\documentclass[12pt,leqno,fleqn]{amsart}
\usepackage{amsmath,amsfonts,amstext,amsthm,amssymb,pxfonts}
\usepackage[colorlinks,citecolor=red,pagebackref,hypertexnames=false]{hyperref}
\usepackage[backrefs]{amsrefs}

\numberwithin{equation}{section}

\allowdisplaybreaks
\swapnumbers

\theoremstyle{plain} 
\newtheorem{lemma}[equation]{Lemma}
\newtheorem{proposition}[equation]{Proposition}
\newtheorem{theorem}[equation]{Theorem}

\newtheorem{smallBallConjecture}[equation]{Small Ball Conjecture}

\newtheorem{sheet}[equation]{The Small Ball Conjecture for the Brownian Sheet}
\newtheorem{product}[equation]{Product Rule in Dimension 2}

\newtheorem{lpi}[equation]{Littlewood Paley Inequalities}

\theoremstyle{definition}
\newtheorem{definition}[equation]{Definition}

\theoremstyle{remark}
\newtheorem*{remark}{Remark}
\newtheorem*{Acknowledgment}{Acknowledgment}

\def\norm#1.#2.{\lVert#1\rVert_{#2}}
\def\Norm#1.#2.{\bigl\lVert#1\bigr\rVert_{#2}}
\def\NOrm#1.#2.{\Bigl\lVert#1\Bigr\rVert_{#2}}
\def\NORm#1.#2.{\biggl\lVert#1\biggr\rVert_{#2}}
\def\NORM#1.#2.{\Biggl\lVert#1\Biggr\rVert_{#2}}

\def\ip#1,#2,{\langle #1,#2\rangle}
\def\Ip#1,#2,{\bigl\langle#1,#2\bigr\rangle}
\def\IP#1,#2,{\Bigl\langle#1,#2\Bigr\rangle}

\def\mid{\,:\,}

\def\abs#1{\lvert#1\rvert}
\def\Abs#1{\bigl\lvert#1\bigr\rvert}
\def\ABs#1{\Bigl\lvert#1\Bigr\rvert}

\def\XXint#1#2#3{{\setbox0=\hbox{$#1{#2#3}{\int}$}
     \vcenter{\hbox{$#2#3$}}\kern-.5\wd0}}

\def\name #1 #2{\operatorname{#1}(#2)}

\begin{document}

\title[Small Ball Inequality in Three Dimensions]
{On the Small  Ball  Inequality in Three Dimensions}

\author[D. Bilyk]{Dmitriy Bilyk}
\author[M.~T.~Lacey]{Michael T. Lacey}

\address{Dmitriy Bilyk\\
School of Mathematics\\
Georgia Institute of Technology\\
Atlanta,  GA 30332 USA}

\email{bilyk@math.gatech.edu}

\address{Michael Lacey\\
School of Mathematics\\
Georgia Institute of Technology\\
Atlanta,  GA 30332 USA}

\email{lacey@math.gatech.edu}

\begin{abstract}
Let $ h_R $ denote an $ L ^{\infty }$ normalized Haar function
adapted to a dyadic rectangle $ R\subset [0,1] ^{3}$.  We show that
there is a positive $ \eta < \frac12$ so that for all integers $ n$,
and coefficients $ \alpha  (R)$ we have
\begin{equation*}
2 ^{-n} \sum _{\abs{ R}=2 ^{-n} } \abs{ \alpha(R) } {}\lesssim{} n
^{1 - \eta   }   \NOrm \sum _{\abs{ R}=2 ^{-n} } \alpha(R) h_R
.\infty . .
\end{equation*}
This is an improvement over the `trivial' estimate by an amount of $
n ^{- \eta }$, while  the Small Ball Conjecture says that  the
inequality should hold with $ \eta =\frac12$. There is a
corresponding lower bound on the $ L ^{\infty }$ norm of the
Discrepancy function of an arbitrary distribution of a finite number
of points in the unit cube in three dimensions.  The prior result,
in dimension $ 3$, is that of J{\'o}zsef Beck \cite{MR1032337}, in
which the improvement over the trivial estimate was logarithmic in $
n$.  We find several simplifications and extensions of Beck's
argument to prove the result above.
\end{abstract}

\maketitle

\section{The Principal Conjecture and the Main Results}
In one dimension, the class of dyadic intervals in the unit interval
is $\mathcal D {} \coloneqq {}\{ [j2^{-k},(j+1)2^{-k})\mid j,k\in
\mathbb N\,,\ 0\le j\le 2 ^{k}-1\} $.
 Each dyadic interval has a left and right half, which are also dyadic.  Define the
 Haar functions
 \begin{equation*}
h_I \coloneqq -\mathbf 1 _{I_{\textup{left}}}+ \mathbf 1
_{I_{\textup{right}}}.
\end{equation*}
 Note that we use an $ L ^{\infty }$ normalization of these functions, which will
 make some formulas seem odd to a reader accustomed to the $L^2$ normalization.

 In dimension $d $, a \emph{dyadic rectangle in the unit cube $ [0,1] ^{d}$}
 is a product of dyadic intervals, thus an
element of $\mathcal D^d $.   A Haar function associated to $R $ is
defined as a product of the Haar functions associated
 with each side of $R $, namely
 \begin{equation*}
 h_{R_1\times\cdots\times R_d }(x_1,\ldots,x_d) {} \coloneqq {}\prod _{j=1}^d h _{R_j}(x_j).
 \end{equation*}
 This is the usual `tensor' definition.

 We will concentrate on rectangles with  fixed volume.
 This is the `hyperbolic'
 assumption, that pervades the subject.
 Our concern is the following Theorem and Conjecture concerning a
 \emph{lower bound} on the $ L ^{\infty }$ norm of sums of hyperbolic Haar functions:

 \begin{theorem}[Talagrand \cite{MR95k:60049},
    Temlyakov \cite{T1}]\label{j.talagrand}
    In dimension $d=2 $, we have
 \begin{equation}  \label{e.talagrand}
2 ^{-n} \sum _{\abs{ R}= 2 ^{-n} } \abs{ \alpha(R) } {}\lesssim{}
\NOrm \sum _{\abs R \ge 2 ^{-n}} \alpha(R) h_R .\infty . .
 \end{equation}
 Here, the sum on the right is taken over all rectangles with area \emph{at least } $ 2
 ^{-n}$.
 \end{theorem}

 \begin{smallBallConjecture} \label{small} For dimension $ d\ge 3$ we have the inequality
 \begin{equation}\label{e.Talagrand}
2 ^{-n} \sum _{\abs{ R}= 2 ^{-n} } \abs{ \alpha(R) } {}\lesssim{} n
^{\frac12(d-2) } \NOrm \sum _{\abs{ R}\ge 2 ^{-n} } \alpha(R) h_R
.\infty . .
\end{equation}
 \end{smallBallConjecture}

 This conjecture is, by one square root of $ n$, better than the trivial estimate available from the
 Cauchy-Schwartz inequality,  see  \S~\ref{s.trivial}.  As well, see that section for
 an explanation as to why the conjecture is sharp.
  The  case of $d=2 $ (with a sum over
  $|R|=2^{-n}$ on the right-hand side) was resolved by Talagrand \cite{MR95k:60049}.
  Temlyakov has given an easier proof of the inequality in its present form \cite{MR96c:41052}, \cite{T1}, which  resonates with the ideas of  Roth
\cite{MR0066435}, Schmidt \cite{MR0319933}, and Hal{\'a}sz
\cite{MR637361}.

Perhaps, it is worthwhile to explain the nomenclature `Small Ball'
at this point. The name comes from the probability theory. Assume
that $X_t: T \rightarrow \mathbb R$ is a canonical Gaussian process
indexed by a set $T$. The \emph{Small Ball Problem} is concerned
with estimates of $\mathbb P ( \sup_{t\in T} |X_t|< \varepsilon )$
as $\varepsilon$ goes to zero, i.e the probability that the random
process takes values in an $L^\infty$ ball of small radius. The
reader is advised to consult a paper by Kuelbs and Li
\cite{MR94j:60078} for a survey of this type of questions. A
particular question of interest to us deals with the Brownian Sheet,
that is, a centered Gaussian process indexed by the points in the
unit cube $[0,1]^d$ and characterized by the covariance relation
$\mathbb E X_s \cdot X_t = \prod_{j=1}^d \min (s_j,t_j)$. The
conjectured form of the aforementioned probability in this case is
the following:

\begin{sheet}
In dimensions $d\ge 2$, for the Brownian Sheet $B$ we have $$-\log
\mathbb P (\norm B.C([0,1]^d). < \varepsilon ) \simeq
{\varepsilon}^{-2} ( \log 1/\varepsilon )^{2d-1}, \quad \varepsilon
\downarrow 0 .$$
\end{sheet}
In dimension $d=2$, this conjecture has been resolved by Talagrand
in the already cited paper \cite{MR95k:60049}, in which he used a
version of \eqref{e.talagrand} for continuous wavelets in place of
Haars to prove the lower bound in the inequality above. In higher
dimensions, the upper bounds are established and the known lower
bounds miss the conjecture by a single power of the logarithm.

Kuelbs and Li \cite{MR94j:60078} have discovered a tight connection
between the Small Ball probabilities and the properties of the
reproducing kernel Hilbert space corresponding to the process, which
in the case of the Brownian Sheet is $WM^2_d$, the Sobolev space of
the functions on $[0,1]^d$ with mixed derivative in $L^2$. In
Approximation Theory, the covering number $N(\varepsilon)$ is
defined as the smallest number of $L^\infty$ balls of radius
$\varepsilon$ needed to cover the unit ball of $WM^2_d$, i.e. the
cardinality of the smallest $\varepsilon$-net, a quantification of
compactness of the unit ball in the uniform metric. The result of
Kuelbs and Li states that

\begin{theorem}
In dimension $d\ge 2$, as $\varepsilon \downarrow 0$ we have
$$-\log
\mathbb P (\norm B.C([0,1]^d). < \varepsilon ) \simeq
{\varepsilon}^{-2} ( \log 1/\varepsilon )^{\beta} \quad \textup{ iff} \quad
\log N(\varepsilon) \simeq {\varepsilon}^{-1} ( \log 1/\varepsilon
)^{\beta/2}.$$
\end{theorem}
This theorem together with Talagrand's work shows that the Small
Ball Conjecture \ref{small} for continuous wavelets implies the
lower bound in the conjectured asymptotics of the covering numbers
$N(\varepsilon)$ (the upper bounds are known). It is also not very
hard to show this implication directly. The Small Ball Conjecture
for the Haar functions implies a lower bound for the covering
numbers of the space $WM^1_d$.  A detailed discussion of the
connections of the Small Ball Conjecture to the Approximation Theory
and other related areas can be found in \cite{MR1005898}, \cite{T2}.

Even though all of the mentioned questions had been completely
resolved in dimension $d=2$, there has been very little progress in
higher dimensions. The main result of the present paper is a partial
resolution of the three dimensional case of the Small Ball
Conjecture. We extend and simplify
 an approach of J.~Beck \cite{MR1032337}, establishing the following
 theorem:

\begin{theorem}\label{t.bl} In dimension $ d=3$,
there is a  positive $ \eta  >0$ for which we have the estimate
\begin{equation}  \label{e.bkl}
2 ^{-n} \sum _{\abs{ R}=2 ^{-n} } \abs{ \alpha(R) } {}\lesssim{} n
^{1 - \eta   }   \NOrm \sum _{\abs{ R}=2 ^{-n} } \alpha(R) h_R
.\infty . .
 \end{equation}
\end{theorem}

 Beck  \cite{MR1032337} established this inequality
 with  $ n ^{-\eta }$ replaced by
 a term logarithmic in $ n$, although Beck  himself
 did not state the result this way, as the principal concern of that paper is on
 the question of Irregularities of Distribution, another area relevant to the Small
 Ball Conjecture.

 In this subject one takes $\mathcal A_N$ to be $N$ points in the $d$-dimensional
 unit cube, and considers
the Discrepancy Function
\begin{equation}  \label{e.discrep}
D_N(x)=\sharp \mathcal A_N \cap [\vec 0,\vec x)-N \abs{ [\vec 0,\vec
x)} .
\end{equation}
Here $[\vec 0,\vec x)=\prod _{j=1}^d [0,x_j)$ is a rectangle with
antipodal corners being $\vec 0$ and $\vec x$. We will typically
suppress the dependence upon the selection of points $\mathcal A_N$.
A set of points will be \emph{well distributed} if this function is
small in some appropriate function space.  Thus, the principal
concern are various lower bounds for the $L^p$ norm of $D_N$. Many
variants of this question are interesting; readers are encouraged to
consult one of the excellent references in this area, e.g.
\cite{MR903025}. The connection\footnote{One expects extremal point
distributions $ \mathcal A_N$ to have about one point in each cube
of volume about $ N ^{-1}$.  Thus the Haar functions adapted to
dyadic rectangles of about this volume are important.} to the Small
Ball Conjecture lies in the `hyperbolic orthogonal function' method
initiated by   Roth \cite{MR0066435} when he proved that for all
dimensions $d\ge2$, $$\norm  D_N.2. \gtrsim (\log
N)^{\frac{d-1}{2}}.$$ Later, Schmidt \cite{MR0319933} has shown that
in dimension 2, the $L^\infty$ norm of the discrepancy function is
much bigger than what the $L^2$ estimate gives us: $$\norm
D_N.\infty. \gtrsim \log N.$$ Notice that, just like in the Small
Ball Conjecture \ref{small}, this beats the $L^2$ bound by one
square root.

Using our method of proof, and well known facts in the literature on
Irregularities of Distribution (\cites{MR903025,MR1032337}), we
obtain
 following theorem:

\begin{theorem}\label{t.beck}
There is a choice of $0<\eta <\tfrac12 $  for which the following
estimate holds for all collections $\mathcal A_N\subset [0,1]^3 $:
\begin{align}  \label{e.beck}
\norm D_N .\infty. {}\gtrsim{} (\log N) ^{1+\eta }\,.
\end{align}
\end{theorem}

Beck's result is as above, with $ (\log N) ^{\eta }$ replaced by  a
doubly logarithmic term in $ N$. There is no further result known to the authors about the
Small Ball Problem, nor the $ L ^{\infty }$ norm of the Discrepancy Function in higher
dimensions.

 Concerning the
 value of $ \eta  $ for which our Theorems hold, it is computable,
 but we do not carry out this step, as the particular $ \eta  $
we would obtain is certainly not  optimal.  Instead, the point of this proof is
that the methods pioneered by J{\'o}zsef Beck are more powerful than originally suspected.
We expect  more efficient organization of the proof, and less \emph{ad hoc } constructions,
will yield quantifiable
and substantive improvements to the results of this paper.\footnote{Additional steps
that one could take to optimize the proof are known to the authors; others are
the subject of speculation.}

 The organization of the proof, at the highest level, and outlined in
  \S~\ref{s.Short}, is that of J{\'o}zsef Beck \cite{MR1032337}.
 At the same time, both the exact construction and subsequent
 details are in many respects easier than in Beck's paper.
 In particular, the
 construction in that section is a Riesz product construction, following the lines
 of \S~\ref{s.talagrand}.  But, the product, with our current understanding,
 must be taken to be `short,' a dictation to us from the
 third dimension:  the `product rule'  \ref{p.product} does not hold
 in dimension three.  This unfortunate, and
critical fact,  forces the
 definition of `strongly distinct' on us. See  Definition~\ref{d.distinct}.
Still, our Riesz product is defined in a way to facilitate the use of
Littlewood Paley inequalities and  conditional expectation arguments,
which is the source of our simplification and strengthening of Beck's
argument.

The principal argument begins in \S~\ref{s.definition}. The earlier
sections of the paper include a brief discussion of prerequisites
for the proof.

\begin{Acknowledgment}
 We have benefited from
several conversations with Mihalis Kolountzakis and Vladimir Temlyakov
on this subject.  A
substantial part of work by the second-named author was done while
in residence at the University of Crete.
\end{Acknowledgment}

 \section{The Trivial Bounds}  \label {s.trivial}

 \emph{Notation.}  The language and notation of probability and
expectation is used throughout.  Thus,
\begin{equation*}
\mathbb E f = \int _{[0,1] ^{d}} f (x) \; dx
\end{equation*}
and $ \mathbb P (A)= \mathbb E \mathbf 1 _{A}$.  This serves to keep
formulas simpler.  As well, certain conditional expectation
arguments are essential to us.  We use the notation
\begin{equation*}
\mathbb P (B\,|\, A)= \mathbb P (A) ^{-1} \mathbb P (A\cap B) \,,
\qquad \mathbb E (B\,|\, A)= \mathbb P (A) ^{-1} \mathbb E (A\cap
B)\,.
\end{equation*}
For a sigma field $ \mathcal F$,
$
\mathbb E (f\,|\, \mathcal F)
$
 is the conditional expectation of $ f$ given $ \mathcal F$.  In all
instances, $ \mathcal F$ will be generated by a finite collection of
atoms $ \mathcal F _{\textup{atoms}}$,  in which case
\begin{equation*}
\mathbb E (f\,|\, \mathcal F)=\sum _{A\in \mathcal F
_{\textup{atoms}}} \mathbb P (A) ^{-1} \mathbb E (f \mathbf 1 _{A})
\cdot \mathbf 1 _{A} \,.
\end{equation*}

We suppress many constants which do not affect the arguments in
essential ways. $ A \lesssim B$ means that there is an absolute
constant so that $ A \le K B$. Thus $ A \lesssim 1$ means that $ A$
is bounded by an absolute constant. And $ A \simeq B$ means $ A
\lesssim B \lesssim A$.

\smallskip

  The inequality (\ref{e.talagrand}) with an extra square root of $ n$ is easy to prove.

 \begin{lemma}\label{l.hyperbolic}
  It is the case that
 \begin{equation*}
 \sum _{\abs R=2 ^{-n}} \abs{\alpha(R)}\cdot \abs R \lesssim
 n ^{\tfrac 12 (d-1)} \NOrm \sum _{\abs R \ge 2 ^{-n}  } \alpha(R) h_R .\infty. \,.
 \end{equation*}
 \end{lemma}

 \begin{proof}  Each point
 $x\in [0,1]^d $, is in at most $  n ^{d-1} $ possible rectangles.  This is the essential
point dictated by the hyperbolic nature of the problem.
 Using this, and the Cauchy--Schwartz inequality, we have
 \begin{align*}
 \sum _{\abs R=2 ^{-n}}
 \abs{ \alpha(R) }\cdot \abs R&{}={}
    \NOrm \sum _{\abs R=2 ^{-n} }
 \abs{ \alpha(R) } {\mathbf 1}_{R} .1.
 \\&{}\lesssim
 n^{\tfrac 12 (d-1)} \NORm\Biggl[ \sum _{\abs R=2 ^{-n}}
 \abs{ \alpha(R) }^2 {\mathbf 1}_{R}\Biggr]^{1/2}  .1.
 \\&
 {}\lesssim  n^{\tfrac 12 (d-1)}\NOrm \sum _{\abs R \ge 2 ^{-n} } \alpha(R) h_R .2.
 \\
 &{}\lesssim  n^{\tfrac 12 (d-1)}\NOrm \sum _{\abs R \ge 2 ^{-n} } \alpha(R) h_R .\infty.
 \end{align*}
 \end{proof}

Let us also see that the Small Ball Conjecture is sharp.
Indeed, we take the $ \alpha (R)$ to be random choices of signs.  It is immediate
that
\begin{equation*}
  2 ^{-n}\sum _{\abs R=2 ^{-n} }
 \abs{ \alpha(R) } \simeq n ^{d-1}\,.
\end{equation*}
On the other hand, for fixed $ x\in [0,1] ^{d}$, by the properties
of Rademacher functions we have
\begin{equation*} \mathbb E \ABs{\sum
_{\abs R=2 ^{-n} } \alpha(R) h_R (x)} \simeq n ^{\tfrac 12
(d-1)}\,.
\end{equation*}
 It is also well known that sums of Rademacher random variables obey a sub--Gaussian
 distributional estimate.  The supremum of such sums admits easily estimated upper bounds.
 In particular, it is enough to test the $ L ^{\infty }$ norm of the sum at a
grid of $ 2 ^{nd}$ points in the unit cube, hence we have
\begin{align*}
\mathbb E \NOrm \sum _{\abs R=2 ^{-n} } \alpha(R) h_R .\infty.
& \lesssim
\sqrt {\log 2 ^{nd}} \cdot
\sup _{x} \mathbb E \ABs{\sum _{\abs R=2 ^{-n}  } \alpha(R) h_R (x)}
 \lesssim n ^{ d/2}\,.
\end{align*}
Comparing these two estimates shows that the Small Ball Conjecture
is sharp. In the trigonometric case, a similar remark has appeared
in \cite{T3}.

 \section{Proof of Talagrand's Theorem} \label{s.talagrand}

 In this section we sketch the proof of V.~Temlyakov \cite{T1} to the stronger inequality
 (\ref{e.talagrand}) in the case of $ d=2$, as this will help understand our construction for $d=3$.
 The line of reasoning is similar to that of Schmidt \cite{MR0319933}.

  The decisive point in two dimensions is that one has a `product rule':

\begin{product}\label{p.product} Let $ R,R'$ be two dyadic rectangles of the
same area.  Then, $h_R \cdot h _{R'}\in \bigl\{ 0\,, \, 1 _{R} \,,\,
\pm h _{R\cap R'} \bigr\}. $ More generally, let $ R_1,R_2,\dotsc,
R_k$ be dyadic rectangles of equal area and distinct lengths in
e.\thinspace g. their first coordinates.  Then $\prod _{j=1} ^{k} h_
{R_j} \in \bigl\{ 0\,,\, \pm h _{R_1\cap \cdots \cap R_k} \bigr\}\,.
$
\end{product}

The fact that this `product rule' fails in higher dimensions is the
most essential complication to the resolution of the Small Ball
Conjecture.

 The proof of (\ref{e.talagrand}) is by duality.  Fix
 \begin{equation*}
H=\sum _{\abs{ R}\ge 2 ^{-n}} \alpha (R) h_R \,.
\end{equation*}
 We will construct a function $\Psi $ with $L^1 $ norm at most $1 $, for which
 the inner product
 \begin{equation}\label{e.HPsi}
 \ip H,\Psi,= 2 ^{-n-1}\sum _{\abs R=2^{-n} } \abs{ \alpha (R)}\,.
 \end{equation}
 This clearly implies Theorem \ref{j.talagrand}.  Moreover, the function $\Psi $ is defined as a
 Riesz product:

 \begin{align*}
 \Psi &  \coloneqq \prod _{s=1}^n (1+\tfrac12 \psi_s )\,,
 \\
 \psi _s& =\sum _{R\,:\,\abs{R_1}=2^{-s}, \abs {R_2}=2^{-n+s} }
\operatorname {sgn} (\alpha (R))  h_R .
 \end{align*}
 Of course $ \Psi $ is non--negative. Moreover, it has $ L^1$ norm one: expanding the
 product, the leading term is $ 1$.  All products of $ \psi _s$ are, by
 Proposition~\ref{p.product}, a sum of Haar functions, hence have mean
 zero. A similar argument implies (\ref{e.HPsi}). The proof is complete.

\section{Littlewood-Paley Theory} \label{s.LP}

In this section we review some basic facts from the Littlewood-Paley
Theory, which will be used repeatedly in subsequent sections. We
state the main inequalities here to make the exposition
self-contained. We also remind the reader that the Haar functions
are normalized to have $ L ^{\infty }$ norm one, so that our
formulas are different from most of our references.

It is important to our applications that we consider the Haar basis
as one for vector valued functions. The vector space should be a
Hilbert space $ \mathcal H$, and by $ L ^{p} _{\mathcal H}$ we mean
the class of measurable functions $ f \mid [0,1] \longrightarrow
\mathcal H$ such  that $ \mathbb E \abs{ f} _{\mathcal H} ^{p} <
\infty $.

The Haar Square Function is
\begin{equation*}
\operatorname S (f) \coloneqq \left[\abs{ \mathbb E f}_{\mathcal H}
^2 + \sum _{I\in \mathcal D}  \frac{{\abs{ \ip f,h_I,} _{\mathcal H}
^2 }}{|I|^2} \, \mathbf 1 _{I}\right] ^{1/2} \,.
\end{equation*}
Here, $\ip f,h_I, =\int_I h_I(x)f(x)dx$ and $\mathbb E f$ should be
understood as Bochner integrals, and we are taking the Hilbert space
norm of those terms that involve $ f$. We shall be applying the
Square Function in the cases when $f$ is a finite linear combination
of Haars, i.e. $f=\sum_{I\in \mathcal I} a_I h_I$, where $\mathcal
I$ is a finite subset of $\mathcal D$ and $(a_I)_{I\in \mathcal I}
\subset \mathcal H$. In this case, $f$ has mean zero and the Square
Function takes the form

\begin{equation*}
\operatorname S (f) = \Bigl[ \sum _{I\in \mathcal I} {\abs{a_I}
_{\mathcal H} ^2 }  \, \mathbf 1 _{I}\Bigr] ^{1/2} \,.
\end{equation*}

Of course we have $ \norm f.2.= \norm \operatorname S (f).2.$ just
due to the fact that $ \{\mathbf 1 _{[0,1]}\}\cup \{h_I \mid I\in
\mathcal D\}$ is an orthogonal basis.

The Littlewood-Paley Inequalities are a  extension of this
equality, to an approximate version that holds on all $ L ^{p}$, $
1<p<\infty $.

\begin{lpi} \label{t.lpi}
For $ 1<p<\infty $ there are absolute constants $ 0<A_p<B_p<\infty $
so that
\begin{equation}\label{e.LLPP1}
\begin{split}
 \norm f.p. &\le B_p \norm \operatorname S (f).p. \,, \qquad 1<p<\infty
 \\
 B_p &\lesssim 1+\sqrt p\, \quad \textup{ for } p\ge2.
\end{split}
\end{equation}
In the reverse direction, we have
\begin{equation} \label{e.LLPP2}
\begin{split}
 A_p \norm \operatorname S (f).p. &\le \norm f .p. \,, \qquad
 1<p<\infty,
 \\
 A_p &\simeq  1+1/\sqrt {p-1}\,.
\end{split}
\end{equation}
\end{lpi}

We stress that these results are delicate.\footnote{To prove our Theorems, we
only need these inequalities with constant $ B_p \lesssim p ^{t}$ for some
fixed power of $ t$.  But, the power of $ t=\frac12$ is the sharp result,
so we use it here.}
Burkholder
\cite{MR976214} has shown that the best constants in the inequality
above for \emph{general martingales} are $ A_p
^{-1}=B_p=\max\{p,q\}-1$.  However, a Haar series is not a general
martingale; it is dyadic, which  forces conditional symmetry. See
\cites{MR1018577,MR1439553,MR800004}.

The constants above are sharp.  To see that $ B_p \simeq \sqrt p$ is
sharp for $ p$ large, just use the Central Limit Theorem for
Rademacher random variables.


\section{Exponential Moments} \label{s.exp}

Let $ \psi \,:\, \mathbb R \longrightarrow \mathbb R $ be a
symmetric convex function with $ \psi (x)=0$ iff $ x=0$.  Define the
Orlicz norm
\begin{equation}\label{e.Orlicz}
\norm f . \psi . \coloneqq \inf \{C>0 \mid \mathbb E \psi (f/C)\le
1\}\,.
\end{equation}
We take the infimum of the empty set to be $ +\infty $, and denote
by $ L ^{\psi  }$ to be the collection of functions for which $
\norm f . \psi .<\infty $. If $ \psi (x)= x ^{p}$, then $ \norm
\cdot .\psi .$ is the usual $ L ^{p}$ norm.

We are especially interested in the class of $ \psi $ given by $
\psi _{\alpha } (x)= \operatorname e ^{\abs{ x} ^{\alpha }}\, ,
\quad \abs{ x} \gtrsim 1 \,$. We will write $ L ^{\psi _\alpha }=
\operatorname {exp} (L ^{\alpha })$. These are the exponential
Orlicz classes. The following equivalence is well known and is based
on Taylor series and Stirling's formula:

\begin{proposition}\label{p.exp}  We have the equivalence of norms
\begin{align*}
\norm f. \operatorname {exp} (L ^{\alpha }).  &\simeq \sup _{p\ge1}
p ^{-1/\alpha } \norm f.p.
 \simeq \sup _{\lambda >0} \lambda ^{\alpha } \lvert\log \mathbb P (
\abs{ f}>\lambda )\rvert\,.
\end{align*}

\end{proposition}

The following distributional estimate holds for hyperbolic sums of
Haar functions:

\begin{theorem}\label{t.distributional} In dimension $ d\ge2$ we have the estimate
\begin{equation}\label{e.distributional}
\NOrm \sum _{\abs{ R}= 2 ^{-n}} \alpha (R) h_R . \operatorname {exp} (L ^{2/ (d-1)}) .
\lesssim
\NOrm \Bigl[ \sum _{\abs{ R}= 2 ^{-n}} \alpha (R) ^2 \mathbf 1 _{R} \Bigr] ^{1/2}
. \infty . \,.
\end{equation}
\end{theorem}

Of principal relevance to us is the three dimensional case, where the estimate
above asserts that the hyperbolic sums are exponentially integrable.

\begin{proof}
The tool is the vector valued Littlewood Paley inequality, with
sharp rate of growth in the constants as $p \to \infty  $, stated in
the previous section. As such the proof is a standard one, see
\cites{MR850744,MR1439553}. We will make use of similar arguments
more than once in this paper.

Applying the one dimensional Littlewood Paley inequality in the coordinate
$ x _1$ we see that
\begin{equation}\label{firstLP}
\NOrm \sum _{\abs{ R}= 2 ^{-n}} \alpha (R) h_R . p . \lesssim \sqrt
p \NOrm \Bigl[  \sum _{r_1=1} ^{n} \Abs{ \sum _{\substack{\abs{ R}=
2 ^{-n}\\ \abs{ R_1}=  2 ^{-r_1} }} \alpha (R) h_R  } ^2 \Bigr]
^{1/2}  .p.
\end{equation}
If we are in dimension $ 2$, note that due to the hyperbolic
assumption, all the rectangles satisfying the conditions of the
summation are disjoint, and thus we have:

\begin{equation} \label{e.=}
\Abs{ \sum _{\substack{\abs{ R}= 2 ^{-n}\\ \abs{ R_1}=  2 ^{-r_1} }}
\alpha (R) h_R  } ^2 =\sum _{\substack{\abs{ R}= 2 ^{-n}\\ \abs{
R_1}=  2 ^{-r_1} }} \abs{ \alpha (R)} ^2   \mathbf 1 _{R} ,
\end{equation}
so our proof is complete in this case.

In the higher dimensional case, the key point is to observe that the
last term can be viewed as an $ \ell ^2 $ space valued function,
that is if we fix all the coordinates except $x_2$ and define an
$\ell_2$-valued function
\begin{equation*}
F(x_2) = \sum_{R_2} \bigg\{ \sum_{\substack{\abs{ R}= 2 ^{-n}\\
\abs{ R_1}=  2 ^{-r_1} }} \alpha (R) \prod_{j\neq 2} h_{R_j} (x_j)
\bigg\} _{r_1=1}^n \, h_{R_2}(x_2),
\end{equation*}
then the expression inside the $L^p$ norm on the right hand side of
(\ref{firstLP}) is exactly $|F|_{\ell^2}$. Thus, the Hilbert space
valued Littlewood Paley inequality applies to the second coordinate,
to give us
\begin{equation*}
\NOrm \sum _{\abs{ R}= 2 ^{-n}} \alpha (R) h_R . p .
\lesssim
 p
\NOrm \Bigl[  \sum _{r_1=1} ^{n} \sum _{r_2=1} ^{n} \Abs{ \sum
_{\substack{\abs{ R}= 2 ^{-n}\\ \abs{ R_j}=  2 ^{-r_j}\,, \ j=1,2 }}
\alpha (R) h_R  }  ^2 \Bigr] ^{1/2}  .p. .
\end{equation*}
Observe that we have a full power of $ p$, due to the two applications
of the Littlewood Paley inequalities.   And if $ d=3$, then analog of (\ref{e.=})
holds, completing the proof in this case.

In the case of dimension $ d\ge 4$ note that we can continue applying
the Littlewood Paley inequalities inductively.  They need only be used $ d-1$
times due to the hyperbolic assumption. Thus, we have the inequality
\begin{equation*}
\NOrm \sum _{\abs{ R}= 2 ^{-n}} \alpha (R) h_R . p.
\lesssim  p ^{ (d-1)/2}
\NOrm \Bigl[ \sum _{\abs{ R}= 2 ^{-n}} \alpha (R) ^2 \mathbf 1 _{R} \Bigr] ^{1/2}
. p .\,, \qquad 2\le p < \infty  \,.
\end{equation*}
The implied constant depends upon dimension; the main point we are interested in
is the rate of of growth of the $ L ^{p}$ norms.
Assuming that the Square Function of the sum is bounded in $ L ^{\infty }$, the $ L ^{p}$
norms can only grow at the rate of $ p ^{(d-1)/2}$, which completes the proof.
\end{proof}

This theorem illustrates a thesis of A.\thinspace Zygmund, which
says that the estimates on product domains are controlled by the
effective number of parameters, which in our hyperbolic setting is
$d-1$. The method of  iteration of the one parameter inequalities,
in the vector valued setting, is a common technique in the subject,
see for instance \cites{MR0252961,MR0290095}. We shall repeatedly
make use of this technique in the present paper.

\section{Definitions and Initial Lemmas for Dimension Three}   \label{s.definition}

As it has been already pointed out, the principal difficulty in
three and higher dimensions is that the product of Haar functions is
not necessarily a Haar function. On this point, we have the
following higher dimensional analogue of the `product rule'
(\ref{p.product}):

\begin{proposition}\label{p.productsofhaars}
Suppose that  $R_1,\ldots,R_k$ are rectangles such that there is no choice of $1\le j<j'\le k$ and no choice of
coordinate $1\le{} t\le d$ for which we have $R _{j,t}=R _{j',t}$.  Then, for a choice of sign $\varepsilon\in
\{\pm1\}$
we have
\begin{equation}  \label{e.productsofhaars}
\prod _{j=1}^k h_R=\varepsilon h_S, \qquad S=\bigcap _{j=1}^k R_k.
\end{equation}
\end{proposition}

\begin{proof}
Expand the product as
\begin{equation*}
\prod _{m=1}^\ell  h_{R_m} (x_1,\dotsc, x_d) = \prod _{m=1}^\ell
  \prod _{t=1} ^{d} h_{R_{m,t}} (x_t).
\end{equation*}
  Our assumption is that for each $ t$, there
is exactly one choice of $ 1\le m_0\le \ell $ such that
$ R _{m_0,t}=S_t$.  And moreover, since the minimum value of $ \abs{ R _{m,t}}$ is obtained
exactly once, for $ m\neq m_0$, we have that
$ h _{R _{m,t}}$ is constant on $ S_t$.  Thus, in the $ t$ coordinate, the product is
\begin{equation*}
  h _{S_t} (x_t)
\prod _{1\le m\neq m_0\le \ell }  h _{R _{m,t}} (S_t) =
\varepsilon_t\, h _{S_t} (x_t),\quad \textup{where $\varepsilon_t\in
\{\pm1\}$} \,.
\end{equation*}
This proves our Lemma.
\end{proof}

\begin{remark}
It is also a useful observation, that the products of Haar functions
have mean zero, if the minimum value of $ \abs{ R _{m,t}}$ is unique
for at least one coordinate $t$.
\end{remark}

Let $\vec r\in \mathbb N^d$ be a partition of $n$, thus $\vec
r=(r_1,r_2,r_3)$, where the  $r_j$ are non negative integers and
$\abs{ \vec r} \coloneqq \sum _t r_t=n$. Denote all such vectors as
$ \mathbb H _n$. (`$ \mathbb H $' for `hyperbolic.') These vectors
will specify the geometry of the rectangles, i.e. we set $\mathcal R
_{\vec r} =\{R\in \mathcal{D}^n :\, |R_j|=2^{-r_j}, \, j=1,2,3 \} $.

We call a function $f$ an $\mathsf r$ function  with parameter $ \vec r$ if
\begin{equation}
\label{e.rfunction} f=\sum_{R\in \mathcal R _{\vec r}}
\varepsilon_R\, h_R\,,\qquad \varepsilon_R\in \{\pm1\}\,.
\end{equation}
We will use $f _{\vec r} $ to denote a generic $\mathsf r$ function.
A fact used without further comment is that $ f _{\vec r} ^2 \equiv 1$.

\begin{definition}\label{d.distinct}  For vectors $ \vec r_j \in \mathbb N ^{3}$,
say that $ \vec r_1,\dotsc,\vec r_J$ are \emph{strongly distinct }
iff for coordinates $ 1\le t\le 3$ the integers $ \{ r _{j,t}\mid
1\le j \le J\}$ are distinct.  The product of strongly distinct $
\mathsf r$ functions is also an $ \mathsf r$ function, which follows
from `the product rule' (\ref{p.productsofhaars}).
\end{definition}

The $\mathsf r$ functions we are interested in are
\begin{equation}\label{e.fr}
f _{\vec r} \coloneqq \sum _{R\in \mathcal R _{\vec r}}
\operatorname {sgn}(\alpha  (R)) \, h_R\,,
\end{equation}
where $ H_n=\sum _{\lvert  R\rvert\ge 2 ^{-n} } \alpha  (R) h_R$.

\section{J{\'o}zsef Beck's  Short Riesz Product}\label{s.Short}

Let us define relevant parameters by
\begin{gather} \label{e.q}
q= a  n ^{\varepsilon}  \,, \qquad   b =\tfrac 16 \, ,
\\
\label{e.rho}
\widetilde \rho=a  q ^{b}  n^{-1}\,, \qquad \rho =  {\sqrt q} n ^{-1}.
\end{gather}
Here, $ a $ are small positive constants, we
use the notation of $ b=\tfrac16$ throughout, so as not to obscure  those
aspects of the argument that that dictate this choice of $ b$.
 $ \widetilde \rho $ is  a `false' $ L^2$
 normalization for the sums we consider, while the larger term $ \rho $ is the
 `true' $ L ^{2}$ normalization.
Our `gain over the trivial estimate' in the Small Ball Conjecture is
$ q ^{b} =n ^{\varepsilon /6}$.  $ 0<\varepsilon <1$ is a small
constant; the exact determination of what we could take $
\varepsilon $ equal to in this proof doesn't seem to be worth
calculating as it surely will not be optimal.

In Beck's paper, the value of $ q=q _{\textup{Beck}}
= \tfrac {\log n} {\log \log n}$ was much smaller than our value of $ q$.
The point of this choice is that $ q
_{\textup{Beck}}^{q _{\textup{Beck}}} \simeq n$, with the term $ q ^{q}$
controlling many of the combinatorial issues concerning the
expansion of the Riesz product.\footnote{Specifically, $ q ^{Cq}$
is a naive bound for the number of admissible graphs, as defined in \S~\ref{s.nsd}.}
With our \emph{substantially  larger value of $ q$}, we need to
introduce additional tools to control the combinatorics.  These
tools are
\begin{itemize}
\item A Riesz product that will permit us to implement various
conditional expectation arguments.
\item Attention to $ L ^{p} $ estimates of various sums, and their growth rates in $ p$.
\item  Systematic use of the Littlewood-Paley inequalities, with the sharp constants in $ p$.
\end{itemize}

Divide the integers $ \{1,2,\dotsc,n\}$ into $ q$ disjoint
increasing intervals $ I_1,\dotsc, I_q$, and let $ \mathbb A _t
\coloneqq \{\vec r\in \mathbb H _n \mid r_1\in I_t\}$.  Let
\begin{equation}
\label{e.G_t}
F_t =  \sum _{\vec r\in \mathbb A _t}  f _{\vec r}\,.
\end{equation}
The Riesz product is now a `short product.'
\begin{equation*}
\Psi \coloneqq \prod _{t=1} ^{q} (1+  \widetilde  \rho F_t) \,.
\end{equation*}
The `false' $ L^2$ normalization implies that the product is, with
high probability, positive, and thus $\norm \Psi.1. \approx \mathbb
E \Psi$, with expectations being typically easier to estimate. This
heuristic is made precise below.

Proposition \ref{p.productsofhaars} suggests that we should
decompose the product $ \Psi $ into
\begin{equation}
\label{e.FG}  \Psi =1+\Psi ^{\textup{sd}}+\Psi ^{\neg  }\,,
\end{equation}
where the two pieces are the `strongly distinct' and `not strongly distinct' pieces.
To be specific, for integers $ 1\le u\le q$,
let
\begin{equation*}
\Psi ^{\textup{sd}}_u
 \coloneqq  \widetilde \rho ^{u}
 \sum _{1\le v_1< \cdots < v_u\le q} \;
 \sideset{}{^{\textup{sd}}}\sum _{\vec r_t\in \mathbb A _{v_t}}  \prod _{t=1} ^{u} f _{\vec r
 _t} \, ,
\end{equation*}
 where $\sideset{}{^{\textup{sd}}}\sum $ is taken to
 be over all $\vec r_t\in  A_{v_t}$ $1\le{} t\le u$ such that:
 \begin{equation}  \label{e.distinct}
 \text{ the vectors  $\{\vec r _{t }\mid 1\le{} t\le u\} $ are strongly
 distinct. }
 \end{equation}
Then define
\begin{equation} \label{e.zCsd}
\Psi ^{\textup{sd}} {} \coloneqq {} \sum _{u=1} ^{q}  \Psi
^{\textup{sd}} _u .
\end{equation}

With this definition, it is clear that we have
\begin{equation} \label{e.gain>trivial}
\begin{split}
\ip H_n, \Psi ^{\textup{sd}}, &=\ip H_n, \Psi ^{\textup{sd}} _{1},
\gtrsim   q ^{b} \cdot n^{- 1}\cdot 2 ^{-n} \sum _{\abs{
R}=2 ^{-n} } \abs{ \alpha(R) }\,,
\\
H_n &= \sum _{\lvert  R\rvert\ge 2 ^{-n} } \alpha (R) h_R\,.
\end{split}\end{equation}
 $ q ^{b}$ is our `gain over the trivial estimate', once we
prove that  $\norm \Psi ^{\textup{sd}}.1. \lesssim 1\,$ (estimate
(\ref{e.sd1}) below). Proving this inequality is the main goal of
the technical estimates of the following Lemma:

\begin{lemma}\label{l.technical} We have these estimates:
\begin{gather}
\label{e.<0} \mathbb P ( \Psi <0) \lesssim \operatorname {exp} (-A
q ^{1/2-b} )\,;
\\  \label{e.expq2b}
\norm \Psi . 2.  \lesssim \operatorname {exp} (a' q ^{2b})\,;
\\
\label{e.E1}
 \mathbb E  \Psi   = 1 \,;
\\
\label{e.CL1}
\norm  \Psi .1. \lesssim 1\,;
\\ \label{e.neq1}
\norm \Psi ^{\neg } .1. \lesssim 1\,;
\\
\label{e.sd1}
\norm \Psi ^{\textup{sd}}.1. \lesssim 1\,.
\end{gather}
Here, $ 0<a'<1$, in (\ref{e.expq2b}), is a small constant, decreasing to zero as $ a $ in (\ref{e.q})
goes to zero; and $ A>1$, in (\ref{e.<0}) is a large constant, tending to infinity as $ a$ in (\ref{e.q})
goes to zero.
\end{lemma}

\begin{proof}
We give the proof of the Lemma, assuming our main inequalities
proved in the subsequent sections.

\smallskip

\emph{Proof of (\ref{e.<0}).} We first note that Theorem
\ref{t.distributional} implies that $\rho F_t$ is in
$\operatorname{exp} (L)$. Then using the distributional estimate of
Proposition \ref{p.exp}, we estimate

\begin{align*} \mathbb P
(\Psi <0)&\le \sum _{t=1} ^{q} \mathbb P ( \widetilde \rho \,  F_t <
-1)
\\
& = \sum _{t=1} ^{q}  \mathbb P (\rho F_t < - a ^{-1} q ^{1/2-b})
\\
& \lesssim  \operatorname {exp} (-c a ^{-1} q ^{1/2-b})\,.
\end{align*}

\smallskip

\emph{Proof of (\ref{e.expq2b}).} The proof of this is detailed
enough and uses the results of subsequent sections, so we postpone
it to Lemma~\ref{l.2b} below.

It is important for our purposes in the proof of the current
Lemma to note that Lemma~\ref{l.2b} proves a uniform estimate, namely
\begin{equation} \label{e.404}
\sup _{V\subset \{1 ,\dotsc, q\}}
\mathbb E \prod _{v \in V}  (1+ \widetilde \rho F_t) ^2 \lesssim
\operatorname {exp} (a' q ^{2b})\,.
\end{equation}

\smallskip

\emph{Proof of (\ref{e.E1}).}
Expand the product in the definition of $ \Psi $.
The leading term is one.  Every other term is a product
\begin{equation*}
\prod _{k\in V} \widetilde \rho \, F_k \,,
\end{equation*}
where $ V$ is a non-empty subset of $ \{1 ,\dotsc, q \}$.  This
product is in turn a linear combination of products of $ \mathsf r$
functions. Among each such product, the maximum in the first
coordinate is unique. This fact tells us that the expectation of
these products of $ \mathsf r$ functions is zero.  So the
expectation of the product above is zero. The proof is complete.

\smallskip

\emph{Proof of (\ref{e.CL1}).} We use the first two estimates of our Lemma.
Observe that
\begin{align*}
\norm \Psi .1. & = \mathbb E \Psi -2 \mathbb E \Psi \mathbf 1 _{\Psi <0}
\\
& \le 1+ 2\mathbb P (\Psi <0) ^{1/2} \norm \Psi .2.
\\
& \lesssim 1+ \operatorname {exp} ( -A q ^{1/2-b}/2+ a' q ^{ 2b})\,.
\end{align*}
We have taken $ b=1/6$ so that $ 1/2-b=2b$.  For sufficiently small
$ a$ in  (\ref{e.q}), we will have $ A \gtrsim a'$.
We see that (\ref{e.CL1}) holds.

In light of the estimate (\ref{e.404}), we see that  the argument above proves
\begin{equation}\label{e.CCLL}
\sup _{V\subset \{1 ,\dotsc, q\}}
\NOrm \prod _{v \in V}  (1+ \widetilde \rho F_t) .1. \lesssim 1\,.
\end{equation}

\smallskip

\emph{Proof of (\ref{e.neq1}).}  The primary facts are (\ref{e.CCLL})
and Theorem~\ref{t.NSD}; we use the notation devised for that Theorem.

Note that the Inclusion-Exclusion principle gives us the identity
\begin{equation*}
\Psi ^{\neg} = \sum _{\substack{V\subset \{1 ,\dotsc, q\}\\  \abs{ V}\ge 2}}
(-1) ^{\abs{ V}+1}
\operatorname {Prod} (\operatorname {NSD} (V))
\cdot \prod _{t\in \{1 ,\dotsc, q\}-V} (1+ \widetilde \rho F_t)\,.
\end{equation*}
We use the triangle inequality, the estimates of Lemma~\ref{l.2b},
H\"older's inequality, with indices $ 1+1/q ^{2b}$ and $ q ^{2b}$,
and the estimate of (\ref{e.NSD}) in the calculation below. Notice
that we have
\begin{align*}
\sup _{V\subset \{1 ,\dotsc, q\}} \NOrm \prod _{v \in V}  (1+
\widetilde \rho F_t) .1+q ^{-2b}. &\le \sup _{V\subset \{1 ,\dotsc,
q\}} \NOrm \prod _{v \in V}  (1+ \widetilde \rho F_t) .1. ^{(1-q
^{-2b})/(1+q ^{-2b})}
\\& \qquad \times
\NOrm \prod _{v \in V}  (1+ \widetilde \rho F_t) .2. ^{q ^{-2b}/(1+q ^{-2b})}
\\
& \lesssim \operatorname{exp} (a'/(1+q^{-2b})) \lesssim 1\,.
\end{align*}

We now estimate
\begin{align*}
\norm \Psi ^{\neg}.1.
&\le
\sum _{\substack{V\subset \{1 ,\dotsc, q\}\\  \abs{ V}\ge 2}} \NOrm
\operatorname {Prod} (\operatorname {NSD} (V))
\cdot \prod _{t\in \{1 ,\dotsc, q\}-V} (1+ \widetilde \rho F_t) .1.
\\
& \le
\sum _{\substack{V\subset \{1 ,\dotsc, q\}\\  \abs{ V}\ge 2}} \norm
\operatorname {Prod} (\operatorname {NSD} (V)) . q ^{2b}. \cdot
\NOrm  \prod _{t\in \{1 ,\dotsc, q\}-V} (1+ \widetilde \rho F_t) .1+ q ^{-2b}.
\\
& \lesssim \sum _{v=2} ^{q} [ q ^{C'}  n ^{-\kappa } ] ^{v}
\lesssim n ^{-\varepsilon '} \lesssim 1\,.
\end{align*}

\smallskip
\emph{Proof of (\ref{e.sd1}).}  This follows from (\ref{e.neq1}) and (\ref{e.CL1})
and the identity $ \Psi =1+\Psi ^{\textup{sd}}+\Psi ^{\neg}$ and the triangle
inequality.

\end{proof}

\section{The Beck Gain in the Simplest Instance} \label{s.beck}

Beck   considered sums of products of $\mathsf r$ functions that are
\emph{not} strongly distinct, and observed that the $ L^2$ norm of
the same are smaller than one would naively expect. This is what we
call the \emph{Beck Gain.} A product of $ \mathsf r$ functions will
not be strongly distinct if the product involves two or more vectors
which agree in  one or more  coordinates. In this section, we study
the sums of products of two $ \mathsf r$ functions which are not
strongly distinct. A later section, \S~\ref{s.nsd}, will study
the general case.  The results of this Section are critical
to the next section, in which we bound the $ L ^2 $ norm of our
Riesz product.

In this section, and again in \S~\ref{s.nsd}, we will use this
notation. For a subset $ \mathbb C \subset \mathbb H _n ^{k}$, let
\begin{equation}\label{e.prod}
\operatorname {Prod} (\mathbb C )
\coloneqq
\sum _{ (\vec r_1 ,\dotsc, \vec r_k) \in \mathbb C }
\prod _{j=1} ^{k} f _{\vec r_j}\,.
\end{equation}
In this section, we are exclusively interested in $ k=2$.

Let
$
\mathbb C (2) \subset  \mathbb H _n ^{2}
$
consist of all pairs of distinct $ \mathsf  r$  vectors
$ \{\vec r_1,   \vec r _{2}\} $ for which
$
 r _{1,2}= r _{2,2}
$.
J.\thinspace Beck calls such terms `coincidences' and we will continue to use that term.
We need norm estimates on
the sums of  products of such $ \mathsf r$ vectors.

\begin{lemma}\label{l.SimpleCoincie}[\textup{\textbf{The  Simplest Instance of the Beck Gain.}}]
We have these
estimates for arbitrary subsets $ \mathbb C \subset \mathbb C (2)$
\begin{align}\label{e.Simple2}
\norm \operatorname {Prod} (\mathbb C ) .p.
& \lesssim   p ^{5/4}  n ^{7/4}\,.
\end{align}
Moreover, if we have $ \mathbb C = \mathbb C (2) \cap \mathbb A _s \times \mathbb A _t$
for some $ 0\le s,t \le q$ we have
\begin{align}\label{e.Simple2q}
\norm \operatorname {Prod} (\mathbb C ) .p. & \lesssim   p ^{3/2}  n
^{3/2}  \,.
\end{align}
\end{lemma}

The second estimate of the Lemma  appears to be sharp, in that the collection $ \mathbb C (2)$
has three free parameters, and the estimates is in terms of $ n ^{3/2}$.
Note that for $ p \simeq n$ we have
\begin{equation*}
\norm \operatorname {Prod} (\mathbb C_2 ) .n.
\simeq
\norm \operatorname {Prod} (\mathbb C_2 ) .\infty . \,.
\end{equation*}
And the latter term can be as big as $ n ^{3}$, which matches the
bound above. Thus we only need to deal with the case $p \lesssim n$.

The proof of the Lemma requires we pass through an intermediary
collection of four tuples of $ \mathsf r$ vectors. Let $ \mathbb B
(4) \subset \mathbb H _n ^{4}$ be four tuples of distinct vectors $
(\vec r,\vec s,\vec t,\vec u)$ for which (i) $ r_2=s_2$ and $
t_2=u_2$; and (ii) in the first and third coordinate  the maximum is
achieved twice.

\begin{proof}
The method of proof is probably best explained by considering first the case
of $ p=2$.  Observe that
\begin{equation*}
\norm \operatorname {Prod} (\mathbb C ).2. ^2 =\mathbb E
\operatorname {Prod} (\mathbb B ) + \mathbb E \operatorname {Prod}
(\widetilde{\mathbb B} )\,,
\end{equation*}
where $ \mathbb B = \mathbb C \times \mathbb C \cap \mathbb B (4) $
and $\widetilde{\mathbb B}$ is a collection of four-tuples in
$\mathbb C \times \mathbb C$ in which some of the vectors completely
coincide. Indeed, the main point is that
\begin{equation*}
\mathbb E f _{\vec r_1} \cdot  f _{\vec r_2} \cdot  f _{\vec r_3} \cdot  f _{\vec r_4}
\neq 0
\end{equation*}
iff the maximum is not unique in each coordinate. But, if the
vectors are distinct, this is the definition of $ \mathbb B (4)$.
Thus the case $ p=2$ follows almost immediately from Lemma~\ref{l.2}
below, since $ \mathbb E \operatorname {Prod} (\widetilde{\mathbb B}
)$ is easy to estimate.

Now, let us consider $ p\ge 4$. Each pair $ (\vec r, \vec s)\in
\mathbb C $ must be distinct in the first and third coordinates.
Therefore, we can apply the Littlewood Paley inequalities in those
coordinates, very much in the same fashion as it was done in the
proof of Theorem \ref{t.distributional}, to estimate
\begin{equation*}
N (p) \coloneqq \norm \operatorname {Prod} (\mathbb C ).p.
\lesssim
 p \NOrm \Bigl[ \sum _{a,b} \ABs{
\sum _{\substack{(\vec r, \vec s)\in \mathbb C \\ \max \{r_1, s_1\}=a
\\ \max \{r_3, s_3\}=b}} f _{\vec r} \cdot f _{\vec s}
} ^2  \Bigr] ^{1/2} .p.\,.
\end{equation*}
Here, we have a full power of $ p$, as we apply the Littlewood Paley inequalities
twice.  Observe that
\begin{equation*}
\sum _{a,b} \ABs{
\sum _{\substack{(\vec r, \vec s)\in \mathbb C \\ \max \{r_1, s_1\}=a
\\ \max \{r_3, s_3\}=b}} f _{\vec r} \cdot f _{\vec s}
} ^2= \sharp \mathbb C + \sum _{i\neq j\in \{1,2,3,4\}}
\operatorname {Prod} (\mathbb C _{i,j})+
 \operatorname {Prod} (\mathbb B )\,.
\end{equation*}
The term $  \sharp \mathbb C$ arises from the diagonal of the
square. The terms $ \mathbb C _{i,j}$ are
\begin{equation*}
\begin{split}
\mathbb C _{i,j}\coloneqq \{ (\vec r_1,\vec r_2,\vec r_3,\vec
r_4)\in \mathbb C \times  \mathbb C \mid  \vec r_i=\vec r_j\,,
\textup{and the other two vectors are distinct}\}.
\end{split}
\end{equation*}
Note that by definition, $\mathbb C _{1,2}=\mathbb C
_{3,4}=\emptyset $, in other cases, the $ \mathbb C _{i,j}$ are of
the same class of objects as $ \mathbb C $. The term $ \mathbb B $
we have already defined.

Then, we can estimate by the triangle inequality,
and the sub-additivity of $ x\mapsto \sqrt x$,
\begin{equation} \label{e.Cinduct}
 p ^{-1}\norm \operatorname {Prod} (\mathbb C ).p.
\lesssim  (\sharp \mathbb C)  ^{1/2} + \sum _{i<j\in \{1,2,3,4\}}
\norm \operatorname {Prod} (\mathbb C _{i,j}) . p/2. ^{1/2} + \norm
\operatorname {Prod} (\mathbb B ). p/2. ^{1/2} \, .
\end{equation}
This inequality is useful for induction.

Let us consider the case of (\ref{e.Simple2q}).  We have already
seen that $ N (2) \lesssim n ^{3/2}  $.   Hence (\ref{e.Cinduct})
implies that for $ p= 2 ^{v+1}$
\begin{equation*}
N (2 ^{v+1}) \lesssim  2 ^{v+1}\bigl\{ n ^{3/2}  + 4 N (2 ^{v})
^{1/2}\bigr\}\,.
\end{equation*}
Clearly, this can be recursively applied, to yield a proof of
(\ref{e.Simple2q}) in the case $p\lesssim n$.
 But the case of $ p\ge n $
is trivial, as the $ L ^{\infty }$ norm of the terms we are estimating
are at most $ n ^{3}$

\end{proof}

\begin{lemma}\label{l.2}    For any subset $ \mathbb B \subset \mathbb B (4)$
\begin{equation}\label{e.2}
\norm \operatorname {Prod} (\mathbb B  ) .p. \lesssim  \sqrt p \, n ^{7/2}\,.
\end{equation}

If we do not consider arbitrary subsets, the estimate improves.  We
have the following
\begin{equation}\label{e.2'}
\norm \operatorname {Prod} (\mathbb B (4)\cap ( \mathbb A _s \times
\mathbb A _t )^{2}  ).p.  \lesssim    p\, n ^{3}\,,
\end{equation}
\end{lemma}

This Lemma, with exponents on $ n$ being $ n ^{7/2}$ appears in
Beck's paper \cite{MR1032337}, in the case of $ p=2$.  The $ L ^{p}$
variants, following from consequences of Littlewood-Paley
inequalities, are important for us.

The first  estimate is recorded, as it is interesting that it
applies to arbitrary subsets of $ \mathbb B (4)$.  We will rely upon
the second estimate.  Pointed out to us by Mihalis Kolountzakis,
this estimate is better for all ranges of $ p \le n$.

\begin{proof}
We discuss (\ref{e.2}). The proof is a case analysis, depending
upon the number of $\{ \vec r, \vec s, \vec t, \vec u\} $ at which
the maximums occur in the first and third coordinates. We proceed
immediately to the cases.

Let $ \mathbb B _2 \subset \mathbb B $ consist of those four--tuples
$ \{ \vec r, \vec s, \vec t, \vec u\}$ for which
\begin{equation*}
r_1=t_1=\max \{r_1,s_1,t_1,u_1\}\,, \qquad r_3=t_3=\max
\{r_3,s_3,t_3,u_3\}\,.
\end{equation*}
This collection is  empty,   for necessarily we must have $
r_2=s_2=t_2=u_2$, but then   $ \vec r=\vec s$, as the parameters of
all vectors is $ n$. This violates the definition of $ \mathbb B $.

\smallskip

Let $ \mathbb B _{3} \subset \mathbb \mathbb B  $ consist of those four--tuples
$ \{ \vec r, \vec s, \vec t, \vec u\}$ for which
\begin{equation*}
r_1=t_1=\max \{r_1,s_1,t_1,u_1\}\,, \qquad r_3=u_3=\max
\{r_3,s_3,t_3,u_3\}\,.
\end{equation*}
That is, the maximal values involve \emph{three} distinct vectors.  These
four vectors can be depicted as
\begin{equation*}
\vec r=\left( \begin{array}{c}
r_1(\Box) \\ r_2  \\ r_3
\end{array}\right) \,,\quad
\vec s= \left( \begin{array}{c} \ast\\
r_2 \\ \Box   \end{array}\right) \,,\quad \vec t= \left(
\begin{array}{c} r_1 \\ t_2 \\ \Box  \end{array}\right) \,,\quad
\vec u= \left( \begin{array}{c} \Box\\ t_2  \\ r_3
\end{array}\right) .
\end{equation*}
A $ \Box$ denotes a parameter which is determined by other choices.
It is essential to note that choices of $ r_2$ and $ r_3$ determine
the value of $ r_1$ (hence the $ \Box$ in the first coordinate for $
\vec r$), and so the vector $ \vec r$.
  The only free parameters
are (say) $ s_1$, denoted by an $ \ast$ above.

But, note that we must then have $ \abs{ \vec s}=s_1+s_2+s_3<n$.  Therefore this
case is empty.

\smallskip

Let $ \mathbb B _{4}$ be those  four-tuples $ \{ \vec r, \vec s,
\vec t, \vec u\}\in \mathbb B$ such that $ s_1=u_1 $ and $ r_3=t_3$.
That is there are \emph{four} vectors involved in the maximums of
the second and third coordinates.  These four vectors can be
represented as
\begin{equation} \label{e.4vec}
\vec r=\left( \begin{array}{c}
\Box \\ r_2  \\r_3\end{array}\right) \,,\quad
\vec s= \left( \begin{array}{c} s_1\\
r_2 \\ \Box   \end{array}\right)  \,,\quad \vec t= \left(
\begin{array}{c} \Box\\ t_2  \\ r_3
\end{array}\right) \,,\quad \vec u= \left(
\begin{array}{c} s_1 \\ t_2 \\ \Box  \end{array}\right) .
\end{equation}

The next argument proves (\ref{e.2}).
Let $ \mathbb B _{4} (a,a',b)$ be those four tuples
$ \{ \vec r,
\vec s, \vec t, \vec u\}\in \mathbb B$ such that
\begin{equation*}
  r_2=s_2=a\,, \quad
  t_2=u_2=a'\, , \quad s_1=u_1=b\,.
\end{equation*}
The point to observe is that
\begin{equation*}
\norm \operatorname {Prod} (\mathbb B _{4} (a,a',b)).p. \le C \sqrt p\,  \sqrt n\,.
\end{equation*}
As there at most $ \lesssim n ^3  $ choices for $ a,a',b$ this will
prove the Lemma.

Indeed, we have not specified $ r_3=t_3$.  Since all vectors are
distinct, we can assume without loss of generality that $ a< a'$
(and thus $r_1>t_1$) and in considering the norm above, we ignore $
\vec s$ and $ \vec u$, as they are completely specified by the datum
$ (a,a',b)$. We apply the Littlewood-Paley inequality in the first
coordinate to the product $ f _{\vec r} \cdot f _{\vec t}$

\begin{equation*}
 \norm f _{\vec r} \cdot f _{\vec t} .p. \lesssim
 \sqrt{p} \NOrm \Bigl[ \sum_{c} \ABs{
\sum _{\substack{\vec r\,, \vec t\,: \\ t_1<r_1=c}} f _{\vec r}
\cdot f _{\vec t}}^2  \Bigr] ^{1/2} .p. = \sqrt{p}\,\sqrt{n}\,,
\end{equation*}
since $r$ and $t$ are completely specified once $r_1$ is fixed. The
proof of (\ref{e.2}) is finished.

\medskip

We turn to the proof of (\ref{e.2'}), arguing similarly. We have
already seen that the only non-empty case is $ \mathbb B _{4}$. Let
$ \mathbb B _{4} (a,a')$ be those four tuples $ \{ \vec r, \vec s,
\vec t, \vec u\}\in \mathbb B_4$ such that
\begin{equation*}
  r_2=s_2=a\,, \quad
  t_2=u_2=a'\,.
\end{equation*}
The point to observe is that
\begin{equation*}
\norm \operatorname {Prod} (\mathbb B _{4} (a,a')).p. \le C  p\,  n\,.
\end{equation*}
As there at most $ \lesssim n  ^2  $ choices for $ a,a'$ this proves the Lemma.

The point is that $ \operatorname {Prod} (\mathbb B _{4} (a,a'))$
almost splits into a product.  Namely, if we define
\begin{align*}
\operatorname {Prod} (\mathbb B _{4,1} (a,a')) & \coloneqq \bigl\{
\{\vec r, \vec t\} \mid  r_2=a,\, t_2=a', \, r_3=t_3 \bigr\}\,,
\\
\operatorname {Prod} (\mathbb B _{4,2} (a,a')) & \coloneqq \bigl\{
\{\vec s, \vec u\} \mid  s_2=a, \, u_2=a', \, s_1=u_1 \bigr\}\,,
\end{align*}
we will have
\begin{equation} \label{e.BAA}
\operatorname {Prod} (\mathbb B _{4} (a,a'))
= \operatorname {Prod}
(\mathbb B _{4,1} (a,a')) \cdot \operatorname {Prod} (\mathbb B
_{4,2} (a,a')) - \operatorname {Prod} (\mathbb M),
\end{equation}
where $\mathbb M \subset (\mathbb B _{4,1} (a,a')) \times (\mathbb B
_{4,2} (a,a'))$ consists of quadruples in which the coincidence
either in the first or the third coordinate is not a maximum in that
coordinate.

We  first prove the  estimate
\begin{equation} \label{e.mk}
\norm
\operatorname {Prod} (\mathbb B _{4,k} (a,a')) .2p.
\lesssim \sqrt p \cdot n ^{1/2}\,, \qquad k=1,2 \,.
\end{equation}
We may assume without loss of generality that $ k=1$, and
$ a>a'$.  The pairs in
$\operatorname {Prod} (\mathbb B _{4,1} (a,a')) $ consist of the two
vectors $ \vec r$ and $ \vec t$ in (\ref{e.4vec}). These two vectors
are parameterized by $ t_1$, say.  Since $a=r_2<a'=t_2 $, and $
r_3=t_3$, the hyperbolic assumption implies $ t_1$ is the maximal
coordinate.  Therefore, the Littlewood-Paley inequality in this
coordinate applies.

Now we deal with the  term $\operatorname {Prod} (\mathbb M)$.
For this, assume that in the first coordinate the maximum is
achieved at $r_1$. This situation is depicted below:

\begin{equation} \label{e.4vecM}
\vec r=\left( \begin{array}{c} \max \\ a  \\r_3\end{array}\right)
\,,\quad
\vec s= \left( \begin{array}{c} s_1\\
a \\ *   \end{array}\right)  \,,\quad \vec t= \left(
\begin{array}{c} * \\ a'  \\ r_3
\end{array}\right) \,,\quad \vec u= \left(
\begin{array}{c} s_1 \\ a' \\ *  \end{array}\right) .
\end{equation}

Notice that in this situation the maximum in the third coordinate
cannot be $r_3=t_3$, for we would then have $s_1+s_2+s_3 < r_1 + r_2
+ r_3 = n$. So, the maximum in this coordinate is $s_3$ or $u_3$.
Also notice, that with $a$ and $a'$ fixed, choosing the values of
$r_1$ and $s_3$ (or $u_3$) completely determines the quadruple of
vectors. Thus we can apply the Littlewood-Paley inequality twice in
the first and the third coordinates, which would yield
\begin{equation} \label{e.MAX}
\Norm \operatorname {Prod} (\mathbb M).p. \lesssim (\sqrt p \sqrt
n)^2 = pn.
\end{equation}

Combining (\ref{e.BAA}), (\ref{e.mk}) and (\ref{e.MAX}),
we see that we have proved
\begin{equation*}
\norm \operatorname {Prod} (\mathbb B _{4} (a,a')) .p.
\lesssim p n\,.
\end{equation*}
The proof is complete.

\end{proof}

There is another corollary to the proof above required at a later
stage of the proof. For an integer $ a$, let $ \mathbb B_a (4)
\subset \mathbb H _n ^{4}$ be four tuples of distinct vectors $
(\vec r,\vec s,\vec t,\vec u)$ for which (i) $ r_2=s_2$ and $
t_2=u_2$; and (ii) in the first  coordinate   we have $ s _{1}=u
_{1}=a$; and (iii) two of the four vectors agree in the third
coordinate.

\begin{lemma}\label{l.Ba}  For any integer $ a$, and
subset $ \mathbb B \subset \mathbb B _{a} (4)$ we have
\begin{equation}\label{e.Ba}
\norm \operatorname {Prod} (\mathbb B ). p. \lesssim p n ^{5/2} \,.
\end{equation}

\end{lemma}

The point of this estimate is that we reduce the number of
parameters of $ \mathbb B (4)$ by one, and gain a full power of $ n$
in the size of the $ L ^{p}$ norm, as compared to the estimate in
(\ref{e.2}).

\begin{proof}
In the proof of Lemma~\ref{l.2}, in the analysis of the terms $
\mathbb B _4$ we used the triangle inequality over the term $
b=s_1=u_1$. Treating this coordinate as fixed, we gain a term $ n
^{-1}$ in the previous proof, hence proving the Lemma above. 

\end{proof}

A further sub-case of the inequality (\ref{e.Simple2}) demands attention.
Using the notation of Lemma~\ref{l.SimpleCoincie}, let
\begin{equation}\label{e.Cadef}
\mathbb C _{2,b} \coloneqq
\{ (\vec r_1, \vec r_2)\in \mathbb C _{2} \mid r _{1,1}=b\}\,,
\qquad 1\le a \le n\,.
\end{equation}
Thus, this collection consists of pairs of distinct vectors, with
a coincidence in the second coordinate, and the first coordinate of $ \vec r_1$
is fixed.
Note that these collections of variables have two free parameters.
At $ L^2$ we find a $ 1/4$ gain over the `naive' estimate.

\begin{lemma}\label{l.Simple2a} For any $ b$
and any subset $ \mathbb C \subset \mathbb C _{2,b}$ we have the estimates
\begin{equation}\label{e.Simple2a}
\norm \operatorname {Prod} (\mathbb C ) .p.
\lesssim   p \cdot  n ^{5/4}\,, \qquad 2 \le p<\infty \,.
\end{equation}

\end{lemma}

\begin{proof}
As in the proof of Lemma~\ref{l.SimpleCoincie}, we begin with the case $ p=2$.
Observer that
\begin{align*}
\norm \operatorname {Prod} (\mathbb C ) .2. ^2
&= \mathbb E \operatorname {Prod} (\mathbb B )\,,
 \end{align*}
 where $ \mathbb B = \mathbb C _{2,b} \times \mathbb C _{2,b}\cap \mathbb B _{b} (4)$,
 with the last collection defined in Lemma~\ref{l.Ba}.
Therefore, the Lemma in this case follows from that Lemma.

More generally, no pair of vectors in $ \mathbb C _{2,b} (2)$ can have a
coincidence in the third coordinate, so we can use the Littlewood Paley
inequalities in that coordinate to estimate
\begin{equation*}
\norm \operatorname {Prod} (\mathbb C ) .p.
\lesssim \sqrt p
\NOrm \Bigl[ \sum _{c} \ABs{ \sum _{\substack{ (\vec r_1, \vec r_2) \in \mathbb C \\
\max \{r _{1,3}, r _{2,3}\} =c} } f _{\vec r_1} \cdot f _{\vec r_2}
} ^2  \Bigr] ^{1/2} .p. .
\end{equation*}

Observe that
\begin{align} \label{e.++}
\sum _{c} \ABs{ \sum _{\substack{ (\vec r_1, \vec r_2) \in \mathbb C \\
\max \{r _{1,3}, r _{2,3}\} =c} }
f _{\vec r_1} \cdot f _{\vec r_2}
} ^2= \sharp \mathbb C
+ \sum _{i<j\in \{1,2,3,4\}} \operatorname {Prod} (\mathbb C _{i,j})
+ \operatorname {Prod} (\mathbb B )\,.
\end{align}
Similar to before, we define the collections $ \mathbb C _{i,j}$ as follows.
\begin{equation*}
\begin{split}
\mathbb C _{i,j}\coloneqq \{ (\vec r_1,\vec r_2,\vec r_3,\vec
r_4)\in \mathbb C \times \mathbb C \mid  \vec r_i=\vec r_j\,,
\textup{and the other two vectors are distinct}\}.
\end{split}
\end{equation*}

In this case, observe that five of these collections are empty, namely
\begin{equation*}
\mathbb C _{1,2}=\mathbb C _{2,3}=\mathbb C _{1,4}=\mathbb C _{2,3}=
\mathbb C _{2,4}=\emptyset\,.
\end{equation*}
The only non-empty collection is $ \mathbb C _{1,3}$. Yet, in $
\mathbb C _{1,3}$, the vectors $ \vec r_2$ and $ \vec r_4$ have a
coincidence in the first coordinate.  Thus,
Lemma~\ref{l.SimpleCoincie} applies to $ \mathbb C _{1,3}$, so that
we have the estimate
\begin{equation} \label{e.13}
\norm \operatorname {Prod} (\mathbb C _{1,3}).p. \lesssim
p ^{5/4} n ^{7/4} \,.
\end{equation}

Let us prove (\ref{e.Simple2a}).
Combining these observations with (\ref{e.++}) and Lemma~\ref{l.Ba} we see that
\begin{align*}
 p ^{-1/2}\norm \operatorname {Prod} (\mathbb C ) .p.
&\lesssim  n +
\norm \operatorname {Prod} (\mathbb C _{1,3}).p/2. ^{1/2}
+ \norm \operatorname {Prod} (\mathbb B ).p/2. ^{1/2}
\\
& \lesssim n +  p ^{5/8} n ^{7/8}+ p ^{1/2} n ^{5/4}\,.
\end{align*}
Concerning the right hand side, note that for $ 2<p< n ^{3}$,
we have $ p ^{5/8} n ^{7/8}< p ^{1/2} n ^{5/4}$. Hence we have proved
\begin{equation*}
\norm \operatorname {Prod} (\mathbb C ) .p.
\lesssim p n ^{5/4}\,, \qquad 1<p < n ^{3}\,.
\end{equation*}
Yet, for $ p \gtrsim n$ the $ L ^{p}$ norm above is comparable to the $ L ^{\infty }$
norm, so we have finished the proof of (\ref{e.Simple2a}).

\end{proof}

\section{The $ L ^{2}$ Norm of the Riesz Product} \label{s.norm}
We now prove a central estimate of the proof.

\begin{lemma}\label{l.2b} The estimate (\ref{e.expq2b}) holds.  Moreover, we have
\begin{equation}\label{e.22bb}
\sup _{V\subset \{1 ,\dotsc, q\}} \mathbb E \prod _{t \in V}  (1+
\widetilde \rho F_t) ^2 \lesssim \operatorname {exp} (a' q ^{2b})\,.
\end{equation}
Here, $ \widetilde \rho $ is as in (\ref{e.rho}), and $ a'$ is a fixed constant times
$ 0<a<1$, the small constant that enters into the definition of $ \widetilde \rho $.
\end{lemma}

\begin{remark}\label{r.condExpect} A conditional expectation argument is essential
to this proof.  This Lemma is also proved in Beck's paper, using a
much more involved argument: his more complicated Riesz product
precludes our simpler line of reasoning.
\end{remark}

\begin{proof}
The supremum over $ V$ will be an immediate consequence of the proof below, and
so we don't address it specifically.

Let us give the initial, essential observation.
We expand
\begin{equation*}
\mathbb E \prod _{t=1} ^{q} (1+ \widetilde \rho F_t) ^2 = \mathbb E
\prod _{t=1} ^{q} (1+ 2\widetilde \rho F_t+ (\widetilde \rho F_t) ^2
)\,.
\end{equation*}
Hold the $ x_2$ and $ x_3$ coordinates fixed, and let $ \mathcal F$ be the sigma
field generated by $ F_1 ,\dotsc, F_{q-1}$.  We have
\begin{equation} \label{e.;p}
\begin{split}
\mathbb E (1+ 2\widetilde \rho F_q+ (\widetilde \rho F_q) ^2
\,\big|\, \mathcal F) &=1+\mathbb E ((\widetilde \rho F_q) ^2
\,\big|\, \mathcal F)
\\
&=1+   a ^2  q ^{2b-1}+ \widetilde \rho ^2  \Gamma _q\,,
\\ \text{where}\,\,
\Gamma _t & \coloneqq \sum _{\substack{ \vec r\neq \vec s \in
\mathbb A _{t}\\ {r_1=s_1} }} f _{\vec r} \cdot f _{\vec s}.
\end{split}
\end{equation}
 Then, we see that
\begin{align} \nonumber
\mathbb E \ \prod _{t=1} ^{q} (1+ 2\widetilde \rho F_t+ (\widetilde
\rho F_t) ^2  ) &= \mathbb E \Bigl\{  \prod _{t=1} ^{q-1} (1+
2\widetilde \rho F_t+ (\widetilde \rho F_t) ^2  )\, \times \mathbb E
(1+ 2\widetilde \rho F_q+ (\widetilde \rho F_q) ^2 \,\big|\,
\mathcal F) \Bigr\}
\\
&\le  \label{e.;;} (1+a ^2 q ^{2b-1}) \mathbb E \prod _{t=1} ^{q-1}
(1+ 2\widetilde \rho F_t+ (\widetilde \rho F_t) ^2  )
\\  \label{e.;;;;}
& \qquad + \Abs{\mathbb E \prod _{t=1} ^{q-1} (1+ 2\widetilde \rho
F_t+ (\widetilde \rho F_t) ^2  ) \cdot
 \widetilde \rho ^{2} \Gamma _q }.
\end{align}
This is the main observation: one should induct on (\ref{e.;;}),
while treating the term in (\ref{e.;;;;}) as an error, as the `Beck
Gain' estimate (\ref{e.Simple2q}) applies to it.

Let us set up notation to implement this line of approach.  Set
\begin{equation*}
N (V;r) \coloneqq \NOrm \prod _{t=1} ^{V} (1+ \widetilde \rho F_t)
.r.  \,, \qquad   V=1 ,\dotsc, q\,.
\end{equation*}
We will use the trivial inequality available from the exponential
moments
\begin{align*}
N (V; 4)&\le \prod _{t=1} ^{V} \norm 1+ \widetilde \rho F_t .4V.
\\
& \le ( 1 +  C q ^{b-1/2} V ) ^{V}
\\
& \le (Cq) ^{Cq}\,.
\end{align*}

This of course is a terrible estimate, but we now use interpolation, noting that
\begin{equation}\label{e.killq^q}
N (V;2(1-1/q) ^{-1} )\le N (V;2) ^{1-1/q} \cdot N (V; 4) ^{1/q}\,.
\end{equation}

We see that (\ref{e.;;}), (\ref{e.;;;;}) and (\ref{e.killq^q})
give us the inequality
\begin{equation}\label{e.==}
\begin{split}
N (V+1; 2) & \le (1+a ^2 q ^{2b-1}) ^{1/2}  N (V; 2)
+ C  \cdot  N (V; 2 (1-1/q) ^{-1} )  \cdot
\norm  \widetilde \rho ^{2} \Gamma _q . q.
\\
& \le
(1+a ^2 q ^{2b-1}) ^{1/2}  N (V; 2)
+ C  N (V;2) ^{1-1/q} \cdot N (V; 4) ^{1/q}
\norm  \widetilde \rho ^{2} \Gamma _q . q.
\\
& \le
(1+a ^2 q ^{2b-1}) ^{1/2}  N (V; 2)
+ C  q ^{C} n ^{-1/2} N (V;2) ^{1-1/q} \,.
\end{split}
\end{equation}
In the last line we have used the inequality (\ref{e.Simple2q}).

Of course we only apply this as long as $ N (V; 2)\ge 1$.  Assuming this is true
for all $ V\ge 1$, we see that
\begin{align*}
N (q;2)&\lesssim  (1+a ^2  q ^{2b-1}+ C q ^{C} n ^{-1/2}) ^{q}
\\
& \lesssim \operatorname e ^{ a' q ^{2b}}\,.
\end{align*}
Here  of course we need $ C q ^{C} n ^{-1/2}\le  a q ^{2b-1} $, which we
certainly have for large $ n$.

\end{proof}

\section{The Beck Gain} \label{s.nsd}

Let us state the main result of this section.
Given $ V\subset \{ 1 ,\dotsc, q\}$ let
\begin{align*}
\operatorname {NSD} (V)
\coloneqq
\Bigl\{ \{\vec r_j \mid j\in V\}\in \operatorname \times _{j\in V} \mathbb A _{j}
\; \big|\;
&\textup{for each $ j\in V$, there is a choice of $ j'\in V- \{j\}$}
\\
& \qquad \textup{and $ \ell =2,3$ so that $ r _{j,\ell }=r _{j',\ell }$}
\Bigr\} \,.
\end{align*}
That is, we take tuples of $ \mathsf r$ vectors, indexed by $ V$, requiring that
each $ \vec r_j$ be in a coincidence.
Such sums admit a favorable estimate  on their $ L^2$ norms.

\begin{theorem}\label{t.NSD}[\textup{\textbf{The Beck Gain.}}]
There are positive constants $ C_0, C_1, C_2, C_3, \kappa    $ for which we have the estimate
\begin{equation}\label{e.NSD}
\rho ^{\abs{ V}} \Norm
\operatorname {Prod} (\operatorname {NSD} (V)) .p.
\lesssim  [C_0 \abs V ^{C_1}  p ^{C_2}  q ^{C_3} n ^{-\kappa  }] ^{\abs{ V}}
\,, \qquad V\subset \{1 ,\dotsc, q\}\,.
\end{equation}
\end{theorem}

\begin{remark}\label{r.gaininV} The novelty in this estimate is that
we find that  (a) the gain can be given in a manner
proportional to $ \lvert  V\rvert $
and (b) the gain also holds in $ L ^{p}$ norms.  In application, $
p \lesssim  q ^{2b}=  q^{1/3} \simeq n ^{\epsilon '}$, so the
polynomial growth in $ p$  and in $ q$ is acceptable to us.\footnote{   Beck
\cite{MR1032337} found a gain in $ L ^{2}$ norm of order $ n
^{-1/4}$, for all  $ V$. Such a small gain of course forces a
much shorter Riesz product.}
\end{remark}

The proof of this Theorem requires a careful analysis of the
variety of ways that a product can fail to be strongly distinct.
That is, we need to understand the variety of ways that coincidences
can arise, and how coincidences can contribute to a smaller  norm.

Following Beck, we will use the language of Graph Theory to describe
these general patterns of coincidences,   although there is
no graph theoretical fact that we need. Rather, the use of this
language is just a convenient way to do some bookkeeping.

The class
of graphs that we are interested in satisfies particular properties.
A \emph{graph} $ G$ is the triple of $ (V (G), E_2, E_3)$, of the
\emph{vertex set} $ V (G)\subset \{1 ,\dotsc, q\}$, and \emph{ edge
sets $ E_2$ and $ E_3$, of color $ 2$ and $ 3$ respectively}.  Edge
sets are are subsets of
\begin{equation*}
E _{j}\subset V (G) \times V (G) - \{ (k,k)\;|\;  k\in V (G)\}\,.
\end{equation*}
Edges are symmetric, thus if $ (v,v')\in E_j$ then necessarily $ (v',v)\in E_j$.

A \emph{clique of color $ j$} is a maximal subset $ Q\subset V (G)$
such that for all $ v\neq v'\in Q$ we have $ (v,v')\in E_j$.  By \emph{maximality},
we mean that no
strictly larger set of vertices  $ Q'\supset Q$ satisfies this condition.

Call a graph $ G$ \emph{admissible} iff
\begin{itemize}
\item  The edges sets, in both colors, decompose into a union of cliques.
\item  Any two cliques  $ Q_2$ in color $ 2$ and clique $ Q_3$ in color $ 3$
can contain at most one common vertex.
\item Every vertex is in at least one clique.
\end{itemize}

 A graph $ G$ is \emph{connected } iff for any two vertices
in the graph, there is a path that connects them.
A \emph{path} in the graph $ G$ is a sequence of vertices  $ v_1 ,\dotsc, v_k$
with an edge of \emph{either color,} spanning adjacent vertices , that is   $ (v _{j}, v _{j+1})
\in E _2 \cup E_3 $.

\subsection*{Reduction to Admissible Graphs}

It is clear that admissible graphs as defined above are naturally
associated to sums of products of $ \mathsf r$ functions. Given
admissible graph $ G$ on vertices $ V$, we set $ X (G)$ to be those
tuples of $ \mathsf r$ vectors
\begin{equation*}
\{\vec r _v \mid v\in V\}\in \prod _{v\in V} \mathbb A _v \, ,
\end{equation*}
so that if $ (v, v')$ is an edge of color $ j$ in $ G$, then $ r _{v,j}= r _{v',j}$.

We will prove the Lemma  below in the following two subsections.
\begin{lemma}\label{l.admissible<} For an admissible graph $ G$ on vertices $ V$ we have
the estimate below for positive, finite constants $ C_0, C_1, C_2, C_3, \kappa $:
\begin{equation}\label{e.equiv}
\rho ^{\abs{ V}} \norm \operatorname {Prod} (X (G)).1.
 \le  [C_0\abs{ V} ^{C_1} p ^{C_2} q ^{C_3}   n ^{-\kappa }] ^{ \abs{ V}}\,,
\qquad 2<p< \infty  \,.
\end{equation}
\end{lemma}

Let us give the proof of Theorem~\ref{t.NSD} assuming this Lemma.
Our tool is the Inclusion-Exclusion Principle, but to apply it we
need additional concepts.

Given two admissible graphs $ G_1, G_2$ on the same vertex set $ V$, let
$ G_1 \wedge G_2$ be the smallest admissible graph which contains all the edges
in $ G_1$ and in $ G_2$.  By smallest, we mean the graph with the fewest number of
edges; and such a graph may not be defined, in which case we take $G_1 \wedge G_2 $
to be undefined.  We recursively define $ G_1\wedge \cdots \wedge G_k \coloneqq
(G_1 \wedge \cdots G _{k-1}) \wedge G_k$.  This wedge product is associative.

Let $ \mathcal G_0$ be the set admissible graphs on $ V$ which are \emph{not}
of the form $ G_1\wedge G_2$ for admissible $ G_1\neq G_2$. These are the `prime' graphs.
(If $ V$ is of cardinality $ 2$ or $ 3$, every graph is prime.)
Now define $ \mathcal G_k$ to be those graphs which are equal to a
wedge product $ G_1 \wedge \cdots \wedge G_k$, with $ G_j\in
\mathcal G_0$, and moreover, $ k$ is the smallest integer for which
this is true. Clearly, we only need to consider $ k\le q$.

Then, by the inclusion-exclusion principle,
\begin{equation}\label{e.INexclude}
\operatorname {Prod} (\operatorname {NSD} (V))
=\sum _{k=0} ^{q} (-1) ^{k}  \sum _{G\in \mathcal G _{k}}
\operatorname {Prod} (X (G))\,.
\end{equation}
The number of admissible graphs
on a set of vertices $ V$ is at most $ 2 ^{\abs{ V}}  \abs{ V}!< 2 ^{\abs{ V}} \abs V ^{\abs
V}$.
So that using (\ref{e.equiv}) clearly implies Theorem~\ref{t.NSD}.

\subsection*{Norm Estimates for Admissible Graphs}

We begin this section with a further reduction to connected admissible graphs.
Let us write $ G\in \operatorname {BG} (C_0, C_1,C_2, C_3 ,\kappa )$ if the
estimates (\ref{e.equiv}) holds.
(`$ \operatorname {BG}$' for `Beck Gain.')
We need to see that all admissible graphs are in $ \operatorname {BG} (C_0, C_1,C_2,
C_3 ,\kappa)$
for non-negative, finite choices of the relevant constants.

\begin{lemma}\label{l.holder} Let $ C_0, C_1,C_2, C_3, \kappa$
be non-negative constants. Suppose that $ G$ is an admissible graph,
and that it can be written as a union  of subgraphs $ G_1 ,\dotsc,
G_k$ on disjoint vertex sets, where all $ G_j \in \operatorname {BG}
(C_0, C_1,C_2, C_3, \kappa)$. Then,
\begin{equation*}
G \in \operatorname {BG} (C_0, C_1, C_2, C_2+C_3, \kappa )\,.
\end{equation*}
\end{lemma}

With this Lemma, we will identify a small class of graphs for which we
can verify the property (\ref{e.equiv}) directly, and then appeal to this
Lemma to deduce Theorem~\ref{t.NSD}.   Accordingly, we modify our notation.
If $ \mathcal G$ is a class of graphs, we write $ \mathcal G\subset \operatorname {BG} (\kappa )$
if there are constants $ C_0, C_1, C_2, C_3 $ such that
$ \mathcal G\subset \operatorname {BG} (C_0, C_1, C_2, C_3, \kappa )$.

\begin{proof}
We then have  by Proposition~\ref{p.products}
\begin{equation*}
\operatorname {Prod} (X (G))
=
\prod _{j=1} ^{k}
\operatorname {Prod} (X (G_j))
\,.
\end{equation*}
Using H\"older's inequality, we can estimate
\begin{align*}
\rho^{|V|}\norm \operatorname {Prod} (X (G)) .p. &\le \prod _{j=1}
^{k} \rho^{|V_j|}\norm \operatorname {Prod} (X (G_j)). k p.
\\
& \le \prod _{j=1} ^{k}  [C_0 (kp) ^{C_1} q ^{C_2} n ^{-\kappa } ] ^{\abs{ V_j}}
\\
& \le [C_0 p ^{C_1} q ^{C_2+C_1} n ^{-\kappa } ] ^{\abs{ V}}\,.
\end{align*}
Here, we use the fact that since the graphs are non-empty, we necessarily have
$ k\le q$.

\end{proof}

\begin{proposition}\label{p.products}
Let $ G_1 ,\dotsc, G_p$ be admissible graphs on pairwise disjoint
vertex sets $ V_1 ,\dotsc, V_p $.  Extend these graphs in the natural way to
a graph $ G$ on the
vertex set $ V=\bigcup V_t$.  Then, we have
\begin{equation*}
\operatorname {Prod} (X(G))
=
\prod _{t=1} ^{p} \operatorname {Prod} (X (G_t))\,.
\end{equation*}
\end{proposition}

\subsection*{Connected Graphs Have the Beck Gain.}

We single out for special consideration the connected  admissible
graphs $ G$ .   Let $ \mathcal G _{\textup{connected}} $ be the
collection of of all admissible  connected
 graphs on $ V\subset \{1,\dotsc, q\}$.

\begin{lemma}\label{l.twoCliques}
We have $ \mathcal G _{\textup{connected}}\subset \operatorname {BG} ( \tfrac1 {15} )$.
\end{lemma}

  We will have
to pay special attention to the case of $ 2$ and $ 3$ vertices . It
is important to observe that the first coordinates are necessarily
distinct, and have the partial order inherited from the vertex set $
V$.   Namely, the vertex set $ V \subset \{1 ,\dotsc, q\}$,
 and $ V$ inherits the order from  the integers.   By the construction
 of our Riesz product, the first coordinates inherit this same order.

\subsubsection*{General Remarks on Littlewood-Paley Inequality.}
These  remarks are essential to our analysis of this lemma, and the Theorem we are proving.
The vertex set $ V$ is a subset of $ \{1 ,\dotsc, q\}$ and it inherits
an order from that set.  Moreover, the tuples of $ \mathsf r$ vectors do as well.
Namely, writing
\begin{equation*}
V= \{v_1< \cdots <v _\ell \},
\end{equation*}
for $ \{\vec r_1 ,\dotsc, \vec r_\ell \} \in X (G)$, we have, by construction,
$ r _{1,1} < \cdots <r _{\ell ,1}$. This since $ r _{m,1} \in I _{v_m}$,
where $ I _{m'}$ is the increasing sequence of intervals of length equal to
$n/q$ that partition $ \{1 ,\dotsc, n\}$.

There is a natural way to apply the Littlewood-Paley inequalities.
For integer $ b _{\ell }\in I _{\ell }$, let $ X (G;   b _{\ell })$
be the tuple of $ \mathsf r$ vectors $ \{\vec r_1 ,\dotsc, \vec r
_{\ell }\}$ such that  $ r _{\ell ,1}=b _{\ell }$.  We have
\begin{equation} \label{e.LPN}
\norm \operatorname {Prod} (X (G)) .p.
\lesssim \sqrt p
\NOrm \Bigl[ \sum _{b _{\ell } \in I _{ v _\ell }}
\abs{ \operatorname {Prod} (X (G;  b _{\ell }))} ^2
\Bigr] ^{1/2} .p. \,.
\end{equation}
It is tempting to continue this procedure, by applying the
Littlewood-Paley inequality again to the vertex $ v _{\ell -1}$.
Yet---and this in an important point---due to the nature of  $
\mathsf r$ functions, this option is blocked to us. The vertex $ v
_{\ell }$ is in at least one clique $ Q $ of, say, color  $ 2$.  We
could choose a value $ c _{Q}$ for that clique, thereby specifying
all coordinates of the vector $ \vec r _{\ell }$.  Set
 $ X (G; b _{\ell } ; c _{Q })$
be the tuple of $ \mathsf r$ vectors $ \{\vec r_1 ,\dotsc, \vec r _{\ell-1 }\}$
such that
\begin{equation*}
\{\vec r_1 ,\dotsc, \vec r _{\ell-1 }\,,\ (b _{\ell }, c_Q, n -  b
_{\ell }-c_Q )\} \in X (G; b _{\ell })\,.
\end{equation*}
Here, $ X (G;  b _{\ell }; c _{q })$
consists of tuples of length $ \ell -1$, since
the vector $ \vec r _{\ell }$ is completely specified.  Thus, we see that
\begin{equation} \label{e.3halfs}
\norm \operatorname {Prod} (X (G)) .p. \lesssim \sqrt p \cdot  n  \sup _{ c_Q}
\NOrm \Bigl[ \sum _{b_ \ell } \operatorname {Prod} (X (G;   b _{\ell
} ; c_Q)) ^2 \Bigr] ^{1/2} .p. \,.
\end{equation}
At this point, the (Hilbert space) Littlewood-Paley inequalities
will again apply.

We will refer to the notation above.  Keep in mind that   $ \vec b$
is for the coordinates specified by a Littlewood-Paley inequality; $
\vec c$ are for the coordinates in a coincidence that we use the
triangle inequality on. We shall return to these themes momentarily.

\begin{proof}[Proof of Lemma~\ref{l.twoCliques}.]
We begin the proof with a discussion of the case of two and three
vertices , which will not be susceptible to the general methods
related to the Littlewood-Paley inequality outlined above.

\subsubsection*{The Case of Two Vertices .}
Notice that if $ G$ consists of only two vertices , the relevant estimate is
(\ref{e.Simple2q}). Namely, we have
\begin{equation*}
\norm \operatorname {Prod} (X (G)) . p. \le C p ^{3/2} n ^{3/2} \,.
\end{equation*}
 Equivalently, $ G\in \operatorname {BG} (C_0, 3/4,0, 1/4)$.

\subsubsection*{The Case of Three Vertices }

The case of $ G\in \mathcal G_{\textup{connected}}$ having three
vertices depends critically on the same phenomena behind the Beck
Gain for graphs on two vertices .  We will deduce this case as a
corollary to the case of two vertices .

There are three distinct sub-cases.  The more delicate of the two
cases is as follows.  The graph is depicted as
\begin{equation}\label{e.badOrder}
\begin{array}{cccccc}
v_1 && v_2 && v_3
\\
\Box  &   & \Box  && \Box
\\
\bullet & = & \bullet  &    &
\\
    &   &  \bullet & = & \bullet
\end{array}
\end{equation}
where $ v_1<v_2<v_3$. (The case of $ v_2<v_1<v_3$ is entirely the same, and we don't
discuss it directly.)

By our general remarks on the Littlewood-Paley inequality, this
inequality applies in the first coordinate, to the vertex $ v_3$.
Using the notation in (\ref{e.LPN}), we have
\begin{equation*}
\norm \operatorname {Prod} (X (G)) .p.
\lesssim \sqrt p
\NOrm \Bigl[ \sum _{b_3 \in I _{v_3 }}
\abs{ \operatorname {Prod} (X (G;   b _{3 }))} ^2
\Bigr] ^{1/2} .p. \,.
\end{equation*}
The vectors $ v_2$ and $ v_3$ have a coincidence in the third coordinate.
Therefore, we specify the value of the coincidence to be $ c_3$ and estimate
\begin{equation} \label{e.badOrder3/2}
\norm \operatorname {Prod} (X (G)) .p.
\lesssim  \sqrt p \cdot   n \cdot  \sup _{c_3}
\NOrm \Bigl[ \sum _{b_3}
\operatorname {Prod} (X (G;  b _{3 }; c_3)) ^2 \Bigr] ^{1/2}   .p. \,.
\end{equation}

Recall that $ X (G;  b_3; c_3)$ consists only of pairs of vectors.
This graph can be depicted as
\begin{equation*}
\begin{array}{cccc}
v_1 && v_2
\\
\Box  &   & \Box
\\
\bullet & = & \bullet
\\
    &   &  c_3
\end{array}
\end{equation*}
But this is the case considered in  (\ref{e.Simple2a}). From that
inequality, we see that we have the estimate
\begin{equation*}
\norm \operatorname {Prod} (X (G;  b _{3 }; c_3)) .p. \lesssim \sqrt
p n ^{5/4} \,.
\end{equation*}
Therefore,
\begin{align*}
\NOrm \Bigl[ \sum _{b_3}
\operatorname {Prod} (X (G;  b _{3 }; c_3)) ^2 \Bigr] ^{1/2}   .p.
& \lesssim \sqrt n  \sup _{b_3} \norm   \abs{ \operatorname {Prod} (X (G;  b _{3 }; c_3))} .p.
\\
& \lesssim \sqrt p \cdot   n ^{7/4}      \,.
\end{align*}
Here we have crudely estimated the $ \ell ^2 $ sum in
(\ref{e.badOrder3/2}).
 Combining the last estimate with (\ref{e.badOrder3/2}), we see
that
\begin{equation} \label{e.badOrder+}
\norm \operatorname {Prod} (X (G)) .p. \lesssim p ^{3/2} n ^{11/4}
\,.
\end{equation}
Recall that the point of comparison is to $\rho^{-3}= n ^{3} q
^{-3/2} $, and the estimate above is smaller by $ n
^{-1/4}$. Thus the class of graphs given by
(\ref{e.badOrder}) are contained in $ \operatorname {BG} (
\tfrac1{12} )$.

\medskip

The other case is when the graph can be depicted by
\begin{equation*}
\begin{array}{cccccc}
v_1 && v_3 && v_2
\\
\Box  &   & \Box  && \Box
\\
\bullet & = & \bullet  &    &
\\
    &   &  \bullet & = & \bullet
\end{array}
\end{equation*}
where $ v_3$, the maximal index is in both cliques.  This case is
much easier, as one application of the Littlewood Paley inequality,
and the triangle inequality will determine the value of both
cliques.   It is very easy to see that this class of graphs is in $
\operatorname {BG} (1/6)$, and the details are omitted. The third
case is even easier -- it involves the graphs which have a clique of
size three in one of the coordinates. Hence the discussion of graphs
on three vertices  is complete.

\subsubsection*{A General Estimate}

We now present a general recursive estimate for the $ L ^{p}$ norm of
$ \operatorname {Prod} (X (G))$, assuming that $ G$ is a connected graph on
at least four vertices.  Write $ V$ as
\begin{equation*}
V= \{v_1< \cdots  < v_\ell \}\,.
\end{equation*}

The estimate is obtained recursively.  Along the way we will construct
two disjoint subsets $ V _{3/2}, V _{1/2} \subset V$.
$ V _{3/2}$ will be the vertices  to which we apply both the Littlewood Paley
and triangle inequalities, thus these vertices  contribute  $ n ^{3/2} q ^{-1/2}$
to our estimate.  $ V _{1/2}$ will be the vertices  to which we apply only
the Littlewood Paley inequality, thus these vertices  contribute $ (n/q) ^{1/2} $
 to our estimate.  Those vertices  not in $ V _{3/2}
 \cup V _{1/2}$ will be those which are determined by earlier steps in the procedure.
 They contribute nothing to our estimate.  In estimating an $ L ^{p}$
 norm, the power of $ p$ is one-half of the number of applications of the
 Littlewood-Paley inequality, namely $ \tfrac 12 \sharp (V _{3/2} \cup V _{1/2})$.

 The purpose of these considerations is to prove the estimate
 \begin{equation} \label{e.purpose}
\norm
\operatorname {Prod} (X (G)) .p.
\le (C \sqrt p) ^{\abs{ V _{3/2}}+ \abs{ V _{1/2}}}
(n/q) ^{ (\abs{ V _{3/2}}+ \abs{ V _{1/2}})/2 } n ^{\abs{ V _{3/2}}} \,.
\end{equation}

Initialize
\begin{equation*}
V _{3/2} \leftarrow \emptyset \,,
\qquad
V _{1/2} \leftarrow \emptyset \,,
\quad
\mathcal Q _{\textup{fixed}} \leftarrow \emptyset \,.
\end{equation*}
The last collection consists of those cliques which are specified by
earlier stages of the argument.

At each stage, we will have an estimate for the form
\begin{equation}\label{e.Estimate}
\begin{split}
\norm
\operatorname {Prod} (X (G)) .p.
&\le (C \sqrt p) ^{\abs{ V _{3/2}}+ \abs{ V _{1/2}}}
n ^{\abs{ V _{3/2}}}
\\
& \qquad \times \sup _{\vec c \in \{ 1 ,\dotsc, n \} ^{\mathcal Q
_{\textup{fixed}}}} \NOrm \Bigl[ \sum _{\vec b \in \{1 ,\dotsc, n \}
^ { V _{3/2}\cup V _{1/2}}} \operatorname {Prod} (X (G;  \vec b;
\vec c)) ^2  \Bigr] ^{1/2} .p. .
\end{split}
\end{equation}

\emph{Base Case of the Recursion.} We update $ V _{3/2} \leftarrow \{v _{\ell }\}$,
since it is the maximal element.  We update $ \mathcal Q _{\textup{fixed}}$
to those cliques which contain $ v _{\ell }$.  Then (\ref{e.Estimate})
is a consequence of (\ref{e.3halfs}).

\emph{Recursive Case.}  At this point, we have the datum $ V
_{3/2}$, $ V _{1/2}$, and $ \mathcal Q _{\textup{fixed}}$. We also
have datum $ \vec b\in  \{1 ,\dotsc, n \} ^ { V _{3/2}\cup V
_{1/2}}$, and $ \vec c \in\{ 1 ,\dotsc, n \} ^{\mathcal Q
_{\textup{fixed}}} $. Notice that this datum can completely specify
some $ \mathsf r$ vectors associated to vertices  not in $ V _{3/2}
\cup V _{1/2}$---think of a vertex that is in two cliques in $
\mathcal Q _{\textup{fixed}}$.

The recursion stops if every  vertex $ v_k$  is determined by this
datum. Otherwise, let $ k$  be the largest integer such that $ \vec
r _{v_k}$ is \emph{not} determined by this datum.   If \emph{no
clique in $ \mathcal Q _{\textup{fixed}}$ contains $ v_k$} update
\begin{equation*}
V _{3/2} \leftarrow V _{3/2} \cup \{v_k\}\,,
\end{equation*}
and update $ \mathcal Q _{\textup{fixed}}$ to include those cliques
which contain $ v_k$.  By application of the Littlewood-Paley
inequality and the triangle inequality, the estimate
(\ref{e.Estimate}) continues to hold for these updated values.

If \emph{some clique in $ \mathcal Q _{\textup{fixed}}$ contains $
v_k$}, then there can be exactly one clique $ Q _{v_k}$ which does,
for otherwise $ \vec r _{v_k}$ would have been completely specified
by these two cliques.  Update
\begin{equation*}
V _{1/2} \leftarrow V _{1/2} \cup \{v_k\}\,,
\end{equation*}
and update $ \mathcal Q _{\textup{fixed}}$ to include  all cliques
which contain $ v_k$.  By application of the Littlewood-Paley
inequality, the estimate (\ref{e.Estimate}) continues to hold for
these updated values.

\smallskip

Once the recursion stops the inequality (\ref{e.Estimate}) holds.  But note that
we necessarily have
\begin{equation*}
\operatorname {Prod} (X (G; \vec b; \vec c)) ^2 \equiv 1\,,
\end{equation*}
as all $ \mathsf r$ vectors are completely determined by $ \vec b$ and $ \vec c$.
Therefore, we have proven (\ref{e.purpose}).

\subsubsection*{The Conclusion of the Proof.}
Since $ V _{3/2} $ and $ V _{1/2}$ are disjoint subsets of $ V$, we have proven
the inequality
\begin{equation}\label{e.Purpose}
\rho ^{\abs{ V}}\norm
\operatorname {Prod} (X (G)) .p.
\le (C \sqrt p) ^{\abs{ V}}
n ^{ \tfrac 32\abs{ V _{3/2}}+ \tfrac 12 \abs{ V _{1/2}} - \abs{ V} }
\,.
\end{equation}
And the remaining analysis concerns the exponent on $ n$ above, namely
we should see that
\begin{equation}\label{e.eta}
\abs{ V} ^{-1} \bigl[ \tfrac 32\abs{ V _{3/2}}+ \tfrac 12 \abs{ V _{1/2}} - \abs{ V} \bigr]
\le  -\tfrac1 {10} \,,
\end{equation}
for a fixed  positive choice of $ \kappa  $, and all connected graphs $ G $
on at least four vertices .  We would conclude that this collection of
graphs  is in $ \operatorname {BG} (\tfrac1 {10} )$.

\smallskip

In order to make the left hand side of (\ref{e.eta}) as large as possible, we should
maximize $ V _{3/2}$.  To continue, we note another formula.
Let $ E (G)$ be the total number of edges in the graph $ G$,
and let $ E (v)$ be the number of edges in $ G$ with one endpoint of the
edge being $ v$.

For $ v\in V _{3/2} \cup V _{1/2}$, let $ F (v)$ be the number of edges
which are specified upon the selection of that vertex in our recursive procedure.
It is clear that we have $ E (v)=F (v)$ if $ v\in V _{3/2}$.  But also,
\begin{equation*}
\sum _{v\in V _{3/2}\cup V _{1/2}} F (v)= E (G)\,.
\end{equation*}
It follows that to maximize the cardinality of $ V _{3/2} $,
those vertices  must be in small cliques.  There are two different classes
of graphs which are extremal with respect to these criteria.

The first  extremal class consists of  graphs $ G$ with all cliques
being of size $ 2$, and the number of cliques is $ \abs{ V}-1$.
For such
graphs, $ \abs{ V _{3/2}} \le \lceil \tfrac 12  \abs{ V}\rceil$, and if
the value is maximal then $ V _{1/2}$ is either $ 0$ if $ \lvert  V\rvert $ is odd,
and $ 1$ if $ \lvert  V\rvert $ is even.  It is straight forward to see that the
maximum of (\ref{e.eta}) occurs at $ \abs{ V}=5$, and is $ -\tfrac1
{10}$. Here, it is vital that we have already discussed the case of
two and three vertices!

The second class are graphs on an even number of vertices, with
half the vertices  in a  clique $ Q$, and each vertex $ v\in Q$ is in
one  clique of size $ 2$.  One can depict such a  graph on six vertices
as
\begin{equation*}
\begin{array}{ccccccccccc}
v_1  && v_2 && v_3 && v_4 && v_5 && v_6
\\
 \ast  &=& \ast &=& \ast
\\
a && b && c && a && b && c
\end{array}
\end{equation*}
The   vertices  are written in increasing order:
$ v_1< v_2< v_3<v_4<v_5<v_6$.  Note that $ v_1,v_2,v_3$
form a single clique of color $ 2$.
There are three additional cliques of size $ 2$, all of color $ 3$. They are
 $ \{v_j,v _{j+3}\}$
for $ j=1,2,3$.
  For such a graph, it is clear that $
\abs{ V _{3/2}}=\tfrac12 \abs{ V}$, and $ \abs{ V _{1/2}}=1$.\footnote{If for example
the maximal vertex $v_6 $ were in the clique of size $ 3$, our algorithm then
predicts a smaller estimate for the graph.}
The
term (\ref{e.eta}) behaves exactly like the first class of extremal  graphs
on an even number of vertices. Our proof is complete.

\end{proof}


\begin{bibsection}
 \begin{biblist}







\bib{MR1032337}{article}{
    author={Beck, J{\'o}zsef},
     title={A two-dimensional van Aardenne-Ehrenfest theorem in
            irregularities of distribution},
   journal={Compositio Math.},
    volume={72},
      date={1989},
    number={3},
     pages={269\ndash 339},
      issn={0010-437X},
    review={MR1032337 (91f:11054)},
}


\bib{MR903025}{book}{
    author={Beck, J{\'o}zsef},
    author={Chen, William W. L.},
     title={Irregularities of distribution},
    series={Cambridge Tracts in Mathematics},
    volume={89},
 publisher={Cambridge University Press},
     place={Cambridge},
      date={1987},
     pages={xiv+294},
      isbn={0-521-30792-9},
    review={MR903025 (88m:11061)},
}




 \bib{MR976214}{article}{
   author={Burkholder, Donald L.},
   title={Sharp inequalities for martingales and stochastic integrals},
   journal={Ast\'erisque},
   number={157-158},
   date={1988},
   pages={75--94},
   issn={0303-1179},
   review={\MR{976214 (90b:60051)}},
}



\bib{MR800004}{article}{
    author={Chang, S.-Y. A.},
    author={Wilson, J. M.},
    author={Wolff, T. H.},
     title={Some weighted norm inequalities concerning the Schr\"odinger
            operators},
   journal={Comment. Math. Helv.},
    volume={60},
      date={1985},
    number={2},
     pages={217\ndash 246},
      issn={0010-2571},
    review={MR800004 (87d:42027)},
}


\bib{MR1439553}{article}{
    author={Fefferman, R.},
    author={Pipher, J.},
     title={Multiparameter operators and sharp weighted inequalities},
   journal={Amer. J. Math.},
    volume={119},
      date={1997},
    number={2},
     pages={337\ndash 369},
      issn={0002-9327},
    review={MR1439553 (98b:42027)},
}

\bib{MR637361}{article}{
   author={Hal{\'a}sz, G.},
   title={On Roth's method in the theory of irregularities of point
   distributions},
   conference={
      title={Recent progress in analytic number theory, Vol. 2},
      address={Durham},
      date={1979},
   },
   book={
      publisher={Academic Press},
      place={London},
   },
   date={1981},
   pages={79--94},
   review={\MR{637361 (83e:10072)}},
}



\bib{MR94j:60078}{article}{
    author={Kuelbs, James},
    author={Li, Wenbo V.},
     title={Metric entropy and the small ball problem for Gaussian measures},
   journal={J. Funct. Anal.},
    volume={116},
      date={1993},
    number={1},
     pages={133\ndash 157},
      issn={0022-1236},
    review={MR 94j:60078},
}





\bib{MR850744}{article}{
    author={Pipher, Jill},
     title={Bounded double square functions},
  language={English, with French summary},
   journal={Ann. Inst. Fourier (Grenoble)},
    volume={36},
      date={1986},
    number={2},
     pages={69\ndash 82},
      issn={0373-0956},
    review={MR850744 (88h:42021)},
}


\bib{MR0066435}{article}{
   author={Roth, K. F.},
   title={On irregularities of distribution},
   journal={Mathematika},
   volume={1},
   date={1954},
   pages={73--79},
   issn={0025-5793},
   review={\MR{0066435 (16,575c)}},
}


\bib{MR0319933}{article}{
   author={Schmidt, Wolfgang M.},
   title={Irregularities of distribution. VII},
   journal={Acta Arith.},
   volume={21},
   date={1972},
   pages={45--50},
   issn={0065-1036},
   review={\MR{0319933 (47 \#8474)}},
}

\bib{MR0252961}{book}{
   author={Stein, Elias M.},
   title={Topics in harmonic analysis related to the Littlewood-Paley
   theory. },
   series={Annals of Mathematics Studies, No. 63},
   publisher={Princeton University Press},
   place={Princeton, N.J.},
   date={1970},
   pages={viii+146},
   review={\MR{0252961 (40 \#6176)}},
}


\bib{MR0290095}{book}{
   author={Stein, Elias M.},
   title={Singular integrals and differentiability properties of functions},
   series={Princeton Mathematical Series, No. 30},
   publisher={Princeton University Press},
   place={Princeton, N.J.},
   date={1970},
   pages={xiv+290},
   review={\MR{0290095 (44 \#7280)}},
}


\bib{MR95k:60049}{article}{
    author={Talagrand, Michel},
     title={The small ball problem for the Brownian sheet},
   journal={Ann. Probab.},
    volume={22},
      date={1994},
    number={3},
     pages={1331\ndash 1354},
      issn={0091-1798},
    review={MR 95k:60049},
}

\bib{MR96c:41052}{article}{
    author={Temlyakov, V. N.},
     title={An inequality for trigonometric polynomials and its application
            for estimating the entropy numbers},
   journal={J. Complexity},
    volume={11},
      date={1995},
    number={2},
     pages={293\ndash 307},
      issn={0885-064X},
    review={MR 96c:41052},
}





\bib{MR1005898}{article}{
   author={Temlyakov, V. N.},
   title={Approximation of functions with bounded mixed derivative},
   journal={Proc. Steklov Inst. Math.},
   date={1989},
   number={1(178)},
   pages={vi+121},
   issn={0081-5438},
   review={\MR{1005898 (90e:00007)}},
}


\bib{T1}{article}{
     author={Temlyakov, V. N.},
     title={Some Inequalities for Multivariate Haar Polynomials},
     journal={East Journal on Approximations},
     volume={1}
     date={1995},
     number={1},
     pages={61 \ndash72}}

\bib{T2}{article}{
     author={Temlyakov, V. N.},
     title={Cubature formulas, discrepancy, and nonlinear approximation},
     journal={J. Complexity},
     number={19},
     date={2003},
     pages={352 \ndash 391},
     }

\bib{T3}{article}{
      author={Temlyakov, V. N.},
      title={An Inequality for Trigonometric Polynomials and its Application for Estimating the Kolmogorov Widths},
      journal={East Journal on Approximations},
      volume={2},
      date={1996},
      pages={253-–262},
      }

\bib{MR1018577}{article}{
   author={Wang, Gang},
   title={Sharp square-function inequalities for conditionally symmetric
   martingales},
   journal={Trans. Amer. Math. Soc.},
   volume={328},
   date={1991},
   number={1},
   pages={393--419},
   issn={0002-9947},
   review={\MR{1018577 (92c:60067)}},
}

\end{biblist}


\end{bibsection}
\end{document}